\documentclass[11pt,reqno]{amsart}
\usepackage{subfiles}
\usepackage{cite}

\usepackage{setspace}
\usepackage[margin=1in]{geometry}
\usepackage[hang,flushmargin,symbol*]{footmisc}
\usepackage{amsmath}
\usepackage{amsthm}
\usepackage{amssymb}
\usepackage{mathtools}
\usepackage{enumitem}
\usepackage{calc}
\usepackage{graphicx}
\usepackage{caption}
\usepackage[labelformat=simple,labelfont={}]{subcaption}
\usepackage{tikz}
\usetikzlibrary{decorations.markings}
\usetikzlibrary{arrows,shapes,positioning}
\usepackage{url}
\usepackage{array}\usepackage{graphicx}
\usepackage{color}

\usepackage{verbatim}

\usepackage[colorinlistoftodos]{todonotes}
\theoremstyle{definition}
\newtheorem{theorem}{Theorem}[section]
\newtheorem{lemma}[theorem]{Lemma}

\newtheorem{corollary}[theorem]{Corollary}
\newtheorem{proposition}[theorem]{Proposition}

\newtheorem{definition}[theorem]{Definition}

\newtheorem{remark}[theorem]{Remark}
\newtheorem{question}[theorem]{Question}

\usepackage{amsfonts}
\usepackage{amsmath}
\usepackage{amssymb}
\usepackage{mathrsfs}
\usepackage{tikz}
\usepackage{tikz-cd}
\usepackage{hyperref}
\usepackage{ mathdots }
\usepackage{dutchcal}
\usepackage{verbatim}
\usepackage{amsmath,amsthm}
\usepackage{extarrows}

\usepackage[mathscr]{euscript}

\newcommand{\Q}{\mathbb{Q}}

\newcommand{\OO}{\mathscr{O}}

\newcommand{\inpt}{\underline{\enspace \thinspace}}

\newcommand{\Hom}{\text{Hom}}
\newcommand{\HHom}{\underline{\text{Hom}}}

\newcommand{\IIsom}{\underline{\text{Isom}}}
\newcommand{\EExt}{\underline{\text{Ext}}}

\newcommand{\Ext}{\text{Ext}}
\newcommand{\Tor}{\text{Tor}}
\newcommand{\hh}[1]{\mathcal{h}_{#1}}

\newcommand{\Amod}{{A \text{-mod}}}
\newcommand{\AAmod}{A' \text{-mod}}

\newcommand{\DDef}{\underline{\text{Def}}}

\newcommand{\coker}{\text{coker}}

\newcommand{\action}{\:\rotatebox[origin=c]{-90}{$\circlearrowright$}\:}

\newcommand{\Split}{{\nabla}}

\tikzset{%
    symbol/.style={%
        ,draw=none
        ,every to/.append style={%
            edge node={node [sloped, allow upside down, auto=false]{$#1$}}}
    }
}

\newcommand{\Arr}{{\text{Arr}}}

\newcommand{\countthm}[1]{
\newcounter{#1}
\setcounter{#1}{\value{theorem}}
}

\newcounter{tempsection}
\newcounter{temptheorem}

\newcommand{\openrefcount}[1]{
\setcounter{tempsection}{\value{section}}
\setcounter{temptheorem}{\value{theorem}}

\setcounter{section}{1}
\setcounter{theorem}{\value{#1}}
}

\newcommand{\closerefcount}{
\setcounter{section}{\value{tempsection}}
\setcounter{theorem}{\value{temptheorem}}
}

\numberwithin{equation}{section}
\numberwithin{paragraph}{section}

\title{\textbf{Deformations of Modules through Butterflies and Gerbes}}
\author{Leo Herr\\
		CU Boulder Math Department}
\date{}

\begin{document}

\maketitle

\begin{abstract}
Classifying obstructions to the problem of finding extensions between two fixed modules goes back at least to L. Illusie's thesis. Our approach, following in the footsteps of J. Wise, is to introduce an analogous Grothendieck Topology on the category $\Amod$ of modules over a fixed ring $A$ in a topos $E$. The problem of finding extensions becomes a banded gerbe and furnishes a cohomology class on the site $\Amod$. We compare our obstruction and that coming from Illusie's work, giving another construction of the exact sequence Illusie used to obtain his obstruction. Our work circumvents the cotangent complex entirely and answers a question posed by llusie. 
\end{abstract}

\section{Introduction}

%Introduction

Consider a topos $E$ and a squarezero extension of sheaves of rings 
\begin{equation}\label{fixalgdef}
0 \rightarrow J \rightarrow A' \rightarrow A \rightarrow 0
\end{equation}
Fix $A$-modules $M$ and $K$, naturally endowed with $A'$-module structures. The central ambition of this paper is to provide another answer to the following question, studied in \cite{illusie1}: 

\begin{question}\label{naivedef}
Is there an extension of $A'$-modules
\begin{equation}\label{moduledef}
\xi : \quad 0 \rightarrow K \rightarrow M' \rightarrow M \rightarrow 0
\end{equation}
and, if so, how many are there? 
\end{question}

We refine Question \ref{naivedef} in two ways. One computes $\Tor_1^{A'}(A, M) = J \otimes_{A'} M$. The boundary map for the long exact sequence of $\Tor$ furnishes an important invariant of the extension $\xi$:  
\begin{equation}\label{inducedmap}
u : J \otimes_{A'} M \rightarrow K
\end{equation}

By sending an extension \eqref{moduledef} to the induced map $u : J \otimes_{A'} M \rightarrow K$, we obtain a morphism
\[\theta : \Ext_{A'}^1(M, K) \rightarrow \Hom_A(J \otimes_{A'} M, K)\]

Extensions may be classified according to their image under $\theta$. In practice, we consider only those extensions which induce a fixed map $u$. 

Given an extension \eqref{moduledef}, we may pull back along a map $N \rightarrow M$ of $A$-modules. 
\[\begin{tikzcd}
\xi|_N :      &0 \arrow{r}      &K \arrow[dashed]{r} \arrow[equals]{d} \arrow[bend left, "0"]{rr}     &M' \times_M N \arrow{r} \arrow{d} \arrow[phantom, very near start, "\ulcorner"]{dr}    &N \arrow{r} \arrow{d}     &0      \\
\xi :     &0 \arrow{r}      &K \arrow{r}     &M' \arrow{r}    &M \arrow{r}     &0
\end{tikzcd}\]

The dashed arrow makes the diagram commute, and the top row is a short exact sequence. This map of extensions is clearly cartesian, forming a fibered category $\EExt_{A'}^1(\inpt, K) \rightarrow \Amod$ over the category of $A$-modules. 

The pullback $\xi|_N$ will map under $\theta$ to the composition $J \otimes_{A'} N \rightarrow J \otimes_{A'} M \overset{u}{\rightarrow} K$. This entails a morphism of fibered categories $\theta : \EExt_{A'}^1(\inpt, K) \rightarrow \Hom_A(J \otimes \inpt, K)$ over $\Amod$, the latter presheaf considered as a fibered category. 

The strict pullback 
\begin{equation}\label{definingdef}
\begin{tikzcd}
\DDef(\inpt, u, K) \arrow{r} \arrow{d} \arrow[phantom, very near start, "\ulcorner"]{dr}      &\EExt_{A'}^1(\inpt, K) \arrow{d}     \\
\{u\} \arrow{r}       &\Hom_A(J \otimes \inpt, K)
\end{tikzcd}
\end{equation}
defines a full fibered subcategory $\DDef(\inpt, u, K) \subseteq \EExt_{A'}^1(\inpt, K)$ of the category of $A'$-module extensions over the category of $A$-modules over $M$, $\Amod/M$. Its sections are extensions $\xi$, with $\theta(\xi)$ a fixed map $u : J \otimes_{A'} M \rightarrow K$. These are referred to as ``deformations.''

Even when it's not possible to construct a deformation of $u$ between $M$ and $K$, it's always possible to find a map $N \rightarrow M$ of $A$-modules and an extension 
\[0 \rightarrow K \rightarrow N' \rightarrow N \rightarrow 0\]
of $A'$-modules which is a deformation of $J \otimes_{A'}N \rightarrow J \otimes_{A'} M \overset{u}{\rightarrow} K$. 

In order to piece together the abundant extensions over $N \rightarrow M$ into one over $M$, we equip $\Amod$ with a topology. 

\begin{definition}[The Topology on $\Amod$]
A family of maps $\{N_i \rightarrow M\}_I$ of modules is deemed covering if, for all finite sets of sections $\Lambda \subseteq M(X)$ over some $X \in E$, there exists a covering $\{U_j \rightarrow X\}_J$ in $E$ so that, for each $j$, there is a single $i$ and a lift of $\Lambda|_{U_j}$ to $N_i(U_j)$. 
\end{definition}

This site is simpler than, but directly analogous to, the site $\OO_Y \text{--} \underline{\text{Alg}}/\OO_X$ of \cite{wisedef1}. The topology is subcanonical. In particular, we write $\hh{K}$ for  the sheaf $N \mapsto \Hom_A(M, K)$. The topology is designed to achieve the next theorem, proved in Section \ref{deformations}. 

\countthm{defisagerbe}
\begin{theorem}\label{defisagerbe}
The fibered category $\DDef(\inpt, u, K) \rightarrow \Amod$ is a gerbe banded by $\hh{K}$. 
\end{theorem}

The fact that $\DDef(\inpt, u, K)$ forms a gerbe answers a few questions for free. 

\begin{theorem}

The class of the $\hh{K}$-gerbe $\DDef(\inpt, u, K)$ in $H^2(\Amod/M, \hh{K})$ obstructs the existence of a deformation $\xi$ with $\theta(\xi) = u$. Provided this class vanishes, the set of such $\xi$ is naturally a torsor under $H^1(\Amod/M, \hh{K})$. The automorphisms of any given extension are in canonical bijection with $H^0(\Amod/M, \hh{K})$. 

\end{theorem}

For any sheaf $F$ on a site $X$, we dogmatically identify $H^1(X, F)$ with $F$-torsors and $H^2(X, F)$ with gerbes banded by $F$ (up to equivalence). This approach provides an instrinsic identification between geometric problems and cohomology classes. There is a choice of sign hidden in this identification -- ours is specified in Definition \ref{torsorgerbesignconvention}. 

In particular, the equivalence class of $\DDef(\inpt, u, K)$ lies naturally in $H^2(\Amod/M, \hh{K})$. The gerbe $\DDef(\inpt, u, K)$ has a global section if and only if it's equivalent to the trivial gerbe $B\hh{K}$, so its equivalence class may be viewed as an obstruction to the existence of an extension. 

If one extension exists, the collection of all extensions is a torsor under $H^1(\Amod/M, \hh{K})$ and the automorphisms of an extension are identified with $\Gamma(\Amod/M, \hh{K})$. This observation answers the ``how many?" of Question \ref{naivedef} in an algebraically refined way. 

The following theorem allows us to compute $H^p(\Amod/M, \hh{K})$. 

\countthm{extiscohom}
\begin{theorem}\label{extiscohom}
The groups of $p$-extensions are all equivalent to cohomology of $\hh{K}$ on the site $\Amod$ : $\Ext^p_A(M, K) \simeq H^p(\Amod/M, \hh{K})$. 
\end{theorem}

This theorem is proved in Section \ref{deformations}. We describe the isomorphism of Theorem \ref{extiscohom} explicitly in the cases $p = 1, 2$ of greatest interest in Propositions \ref{explicitp1} and \ref{explicitsplittingsp2}. As a result of this description, we can identify which 2-extension corresponds to our gerbe $\DDef(\inpt, u, K)$ in Section \ref{extandcohom}: 

\countthm{obstcomparison}
\begin{theorem}\label{obstcomparison}
The diagram 
\begin{equation}\label{obstcomparisondiagram}
\begin{tikzcd}
\Hom_A(J \otimes_{A'}M, K) \arrow{dr}[swap]{\DDef} \arrow{r}{\smile \omega}      &\Ext_A^2(M, K) \arrow{d}{\Split}     \\
        &H^2(\Amod/M, \hh{K})
\end{tikzcd}
\end{equation}
anti-commutes. 
\end{theorem}

The arrow $\Split$ is the isomorphism of Theorem \ref{extiscohom}, $\smile \omega$ sends a map $f : J \otimes_{A'}M \rightarrow K$ to the pushout of a given 2-extension $\omega$ along $f$, and $\DDef$ sends a morphism $f : J \otimes_{A'}M \rightarrow K$ to the $\hh{K}$-gerbe $\DDef(\inpt, f, K)$ over $\Amod/M$. Anti-commutativity signifies that $\DDef(\inpt, f, K)$ and $\Split(f \smile \omega)$ represent additive-inverse cohomology classes. In other words, our obstruction and Illusie's are inverses.

The classification of deformations found in \cite{illusie1} produces the complex
\begin{equation}\label{illusieexactseq}
0 \rightarrow \Ext_A^1(M, K) \rightarrow \Ext_{A'}^1(M, K) \overset{\theta}{\rightarrow} \Hom_A(J \otimes_{A'}M, K) \overset{\smile \omega}{\rightarrow} \Ext_A^2(M, K)
\end{equation}

\begin{lemma}
The sequence of maps \eqref{illusieexactseq} is an exact sequence. 
\end{lemma}

\begin{proof}

Choose an extension
\[\xi : \quad 0 \rightarrow K \rightarrow M' \rightarrow M \rightarrow 0 \quad \in \Ext_{A'}^1(M, K)\]

The action of $J$ on $M'$ factors as $J \otimes M' \twoheadrightarrow J \otimes M \overset{\theta(\xi)}{\rightarrow} K \hookrightarrow M'$ by the definition of $\theta$. Since the first map is surjective and the last is injective, the composite is zero precisely when $\theta(\xi)$ is. Observe that $J$ annihilates $M'$ if and only if $M'$ is an $A$-module if and only if $\xi \in \Ext_A^1(M, K)$. This proves exactness at the domain of $\theta$. 

By Theorem \ref{obstcomparison}, $\DDef(\inpt, u, K)$ and $\Split(u \smile \omega)$ are inverse cohomology classes. One gerbe has a section when the other does. For the gerbe $\DDef(\inpt, u, K)$ to have a section, $u$ must be in the image under $\theta$ of some extension. For $\Split(u \smile \omega)$, this means that the 2-extension $u \smile \omega$ is equivalent to zero in $\Ext_A^2(M, K)$. 

\end{proof}

Exactness entails that the pushout $f \smile \omega$ is equivalent to the zero 2-extension precisely when $f$ is in the image of $\theta$. Under this light, Theorem \ref{obstcomparison} says Illusie's obstruction $f \smile \omega$ is identified with the inverse of our $\DDef(\inpt, f, K)$ under the isomorphism $\Split$. This answers a generalization of Question 3.1.10 in \cite{illusie1}. 

The exact sequence \eqref{illusieexactseq} originates in the transitivity triangle for the graded cotangent complex produced in \cite{illusie1}. Our concrete descriptions of the maps augment those found in \cite[Tag 08L8]{sta}. We can also construct the sequence without reference to the cotangent complex as follows.

Restrict scalars along the map $A' \rightarrow A$ to get a fully faithful embedding $r : \Amod/M \rightarrow \AAmod/M$. This is how we consider $M$ and $K$ as $A'$-modules, and we often continue to suppress the notation $r$. The functor $r$ is cover-preserving and left exact, yielding a morphism of sites
\[\pi : \AAmod/M \rightarrow \Amod/M\]

To avoid ambiguity, we write cohomology on $\Amod$ as $H^p(A/M, \hh{K})$ and that on $\AAmod$ as $H^p(A'/M, \hh{K})$ (and similarly for global sections). The equality $\Gamma(A'/M, \hh{K}) = \Gamma(A/M, \pi_* \hh{K})$ witnesses that the two global section maps to $(Sets)$ commute. The Grothendieck-Leray Spectral Sequence 
\begin{equation}
E_2^{p, q} : \quad H^p(A/M, R^q \pi_* \hh{K}) \Rightarrow H^{p+q}(A'/M, \hh{K})
\end{equation}
yields a 5-term exact sequence. The concern of Section \ref{illusiecomparison} is the next theorem.  

\countthm{illusie5term}
\begin{theorem}\label{illusie5term}
Illusie's exact sequence \eqref{illusieexactseq} and the 5-term exact sequence from the Grothendieck-Leray spectral sequence are isomorphic. The diagram with Illusie's exact sequence on top and the 5-term exact sequence on the bottom commutes: 
\begin{equation}\label{illusie5termiso}
\begin{tikzcd}
0 \arrow{r}       &\Ext_A^1(M, K) \arrow{r} \arrow{d}{\sim}     &\Ext_{A'}^1(M, K) \arrow{r} \arrow{d}{\sim}    &\Hom_A(J \otimes_{A'}M, K) \arrow{r} \arrow[dashed, "\sim"]{d}     &\Ext_A^2(M, K) \arrow{d}{\sim}             \\
0 \arrow{r}       &H^1(A/M, \hh{K}) \arrow{r}      &H^1(A'/M, \hh{K}) \arrow{r}     &H^0(A/M, R^1\pi_* \hh{K}) \arrow{r}       &H^2(A/M, \hh{K}) 
\end{tikzcd}
\end{equation}
\end{theorem}

The solid vertical arrows of \eqref{illusie5termiso} are the isomorphisms of Corollary \ref{extiscohom}, except the last one has a minus sign. Once we show the diagram is natural in $M$ and $K$ in Lemma \ref{illusie5termisnatural}, we obtain the dashed arrow by sheafifying in $M$. Immediately from this identification, we may extend Illusie's exact sequence to the right via $\Ext_A^2(M, K) \rightarrow \Ext_{A'}^2(M, K)$. 

Our present work most heavily relies on the paper \cite{wisedef1}. However, there is an error in the proof which we will correct in later work (see Remark \ref{counterexamplewisedef}). This paper is logically independent of \cite{wisedef1} and none of the present article depends on the mistaken assertions therein.

In the body of the paper, we omit the subscript $\otimes_{A'}$ and write $\pi_* \EExt_{A'}^1(\inpt, K)$ for the strict pullback of $\EExt_{A'}^1(\inpt, K) \rightarrow \AAmod$ along $\pi^*$ as in \cite[04WA]{sta} (for brevity and clarity, respectively).

The material of this paper will likely belong to the author's thesis from CU Boulder under the supervision of Jonathan Wise. The author would like to thank J. Wise for his patience and insight as well as the enormity of his contribution to the present work. Almost all the results here began with him and were developed together in constant communication.

\section{The Topology on A-mod}

%The topology on $\Amod$ with some remarks

In this section, we prove that $\EExt_A^1(\inpt, K) \rightarrow \Amod$ is a stack. We collect a number of convenient properties of $\Amod$ along the way.

For an object $S \in E$, write $A^S$ for the sheafification of the presheaf $U \mapsto \bigoplus_{\Gamma(U, S)} \Gamma(U, A)$. It deserves the title ``free module'' via universal property. 

Suppose $j : U \rightarrow *_E$ is a map to the final object, and $\Lambda \subseteq \Gamma(U, M)$ is a finite subset. The constant sheaf $\underline{\Lambda}$ on $U$ has an adjoint map $\underline{\Lambda} \rightarrow M|_U$, and another adjunction furnishes $j_! \underline{\Lambda} \rightarrow M$ in $E$. 

One particularly useful tautological cover by free modules is $\{A^{j_!\underline{\Lambda}} \rightarrow M\}$, ranging over all such finite subsets of sections. Another is the single element cover $A^M \rightarrow M$. A defining characteristic of our topology is the availability of such covers by free modules.

\begin{lemma}\label{subcanonamod}
The topology on $\Amod$ is subcanonical. 
\end{lemma}

\begin{proof}

Let $\{M_i \rightarrow M\}_I$ be a cover. We check by hand that 
\begin{equation}\label{babycech}
\bigoplus_{I \times I} (M_i \times_M M_j) \rightarrow \bigoplus_I M_i \rightarrow M \rightarrow 0
\end{equation}
is exact. The leftmost arrow is the difference of the two projections. 

For $\hh{K}$ to be a sheaf, the complex obtained by applying $\hh{K}$ to this one must be exact. Left exactness of $\hh{K}$ will give us the result. It's clear that the sequence is a complex and that $\bigoplus_I M_i \rightarrow M$ is surjective. 

First assume the cover consists of a single element, $\{T \rightarrow M\}$. Let $S$ be the kernel, fitting into a short exact sequence
\[0 \rightarrow S \rightarrow T \rightarrow M \rightarrow 0\]

Remark that $T \oplus S \rightarrow T \times_M T$ sending local sections $(t, s) \mapsto (t + s, t)$ is an isomorphism. The composite of this isomorphism with the map $T \times_M T \rightarrow T$ of \eqref{babycech} is projection onto $S$ and inclusion. Then \eqref{babycech} takes the form
\[T \oplus S \rightarrow T \rightarrow M \rightarrow 0,\]
which is exact. 

Now return to the general case of an arbitrary cover. In order to verify \eqref{babycech} is exact, we may freely localize in $E$. We argue that any local section of $\bigoplus_I M_i$ is cohomologous to a section of a single $M_0$ (among the $M_i$) locally in $E$, reducing the verification of exactness to the special case considered above. 

After localization in $E$, all sections of $\bigoplus_I M_i$ are represented by finite sums of sections from various $M_i$. Choose $\sum \limits_{k=1}^n m_{i_k} \in \Gamma(U, \bigoplus_I M_i)$ and consider the images $\overline{m}_{i_k}$ of $m_{i_k}$ in $M$. Localize in $E$ again and use the covering condition to lift the finite set of sections $\{\overline{m}_{i_k}\}_{k=1}^n \subseteq \Gamma(U, M)$ to some single $M_0$ among the $M_i$. Let $m'_{i_k}$ be a chosen preimage in $M_0$ of $\overline{m}_{i_k}$. 

Consider the section $\sum \limits_{k=1}^n (m'_{i_k}, m_{i_k})$ of $\bigoplus \limits_{I \times I} (M_i \times_M M_j)$. The second projection maps this section to the one we started with; the first yields a sum of elements of $M_0$. Hence our original section is cohomologous to one in $M_0$.

\end{proof}

Our next goal is to show the topology makes the fibered category $\EExt_A^1(\inpt, K) \rightarrow \Amod$ into a stack. Recall that extensions up to isomorphism form a group, with identity given by the trivial extension 
\[\underline{0} :  \quad 0 \rightarrow K \rightarrow K \oplus M \rightarrow M \rightarrow 0\]
The trivial extension is isomorphic to any extension whose epimorphism admits a section. Addition (the ``Baer Sum'') of two extensions 
\begin{equation}\label{extensionnames}\xi : \quad 0 \rightarrow K \rightarrow M' \rightarrow M \rightarrow 0
\quad \quad \eta : \quad 0 \rightarrow K \rightarrow M'' \rightarrow M \rightarrow 0\end{equation}
is defined by pulling back and pushing out the product of the two extensions along the maps in the diagram: 
\[\begin{tikzcd}
        &       &       &M \arrow{d}{(id, id)}      &       \\
0 \arrow{r}       &K \oplus K \arrow{r} \arrow{d}{id + id}      &M' \oplus M'' \arrow{r}     &M \oplus M \arrow{r}     &0      \\
        &K      &       &       &
\end{tikzcd}\]

In other words, the group law is defined by biadditivity and functoriality: 
\[\Ext_A^1(M, K) \times \Ext_A^1(M, K) \rightarrow \Ext_A^1(M \oplus M, K \oplus K) \rightarrow \Ext_A^1(M, K)\]
We fix the notation $\xi$ and $\eta$ for the extensions above throughout this section.

\begin{remark}\label{basicprops}
We collect a few basic properties of the topology on $\Amod$. 
\begin{itemize}
\item For $N, N' \in \Amod/M$, the presheaf $\HHom_A(N, N')$ sending $P \mapsto \Hom_A(P \times_M N, P \times_M N')$ and the subpresheaf of isomorphisms $\IIsom_A(N, N')$ are both sheaves. 

\item Extensions $\xi$ as in \eqref{extensionnames} are locally isomorphic to the trivial extension over $\Amod$. 

\item Given two families of maps $\{N_i \rightarrow M\}$ and $\{P_j \rightarrow M\}$, if the latter is covering and refines the former via maps $\{P_j \rightarrow N_{i_j}\}$ over $M$, then the former is also covering. 
\end{itemize}

The first point follows formally from the subcanonicity of the topology. The second is shown by pullback along $M' \rightarrow M$ and the third follows from the definition of the topology. 

\end{remark}

Recall the trivial gerbe $B\hh{K} \rightarrow \Amod$ whose sections over some $M$ are $\hh{K}$-torsors: sheaves $P$ on $\Amod/M$ which carry a free and transitive action of $\hh{K}|_M$ and locally admit sections. We often write $\hh{M'|M}$ for the sheaf some $M' \in \Amod/M$ represents to emphasize the structure map $M' \rightarrow M$, as opposed to the sheaf $\hh{M'}$ on $\Amod$.

\begin{definition}\label{definerho}
The functor $\rho : \EExt_A^1(\inpt, K) \rightarrow B\hh{K}$ between fibered categories sends an extension $\xi$ to $\hh{M'|M}$. This sheaf becomes a $\hh{K}$-torsor via addition 
\[K \times M' \rightarrow M' \times M' \overset{+}{\rightarrow} M'\]
and the Yoneda Embedding. Morphisms of extensions induce morphisms of representable sheaves which are $\hh{K}$-equivariant. By Remark \ref{basicprops}, the sheaves $\hh{M'|M}$ are isomorphic to their structure group $\hh{K \oplus M} \simeq \hh{K}|_M$ after pullback along a cover in $\Amod/M$. The Yoneda Lemma verifies that $\rho$ is fully faithful. 
\end{definition}

We will show $\rho$ is an equivalence. Our construction relies on the free module functor
\[\lambda : E \rightarrow \Amod\] 
sending $S \mapsto A^S$.

\begin{lemma}\label{finlimfreecovers}
The free module functor $\lambda$ sends fiber products to covers. That is, the natural map 
\[A^{S \times_R T} \rightarrow A^S \times_{A^R} A^T\]
is covering, for $S, R, T \in E$. 
\end{lemma}

\begin{proof}

Since the family we wish to show is a cover consists of a single map, it suffices to show it's a cover in $E$ instead of $\Amod$. Choose a section $\alpha \in \Gamma(U, A^S \times_{A^R} A^T)$; we wish to find a lift of $\alpha$ to $A^{S \times_R T}$ locally in $E$. Locally, we may assume $\alpha = (\sum x_k s_k, \sum y_k t_k)$ is a pair of finite sums with $x_k, y_k \in A(U)$, $s_k \in S(U)$, $t_k \in T(U)$ with the same image in $A^R(U)$. 

Fix $r \in R(U)$ and suppose the $s_k$, $t_k$ mapping to $r$ are numbered $\{s_1, \cdots, s_n\}$, $\{t_1, \cdots, t_n\}$. In order for the two sums to have the same image, we must have 
\[\sum \limits_{k = 1}^n x_k = \sum \limits_{k = 1}^n y_k \in A(U)\]
Define $z$ to be the value of either sum. 

Consider the section $\beta_r$ of $A^{S \times_R T}$ given by the sum of $(s_i, t_j)$ with coefficients
\[ \begin{cases} 
    z   & i = j = 1 \\
    x_i + y_j   & i = j \neq 1  \\
    -y_i    & j = 1 \neq i  \\
    -x_j    & i = 1 \neq j  \\
    0       & \text{otherwise}
   \end{cases}
\]

Writing the coefficients as a matrix yields 
\[ \left( \begin{array}{cccc|c}
z & -y_2 & \cdots   &-y_n    &t_1\\
-x_2 & x_2 + y_2 &  &0    &t_2\\
\vdots &  & \ddots &    &\vdots\\
-x_n    &0   &   &x_n + y_n    &t_n\\
\hline s_1&s_2&\cdots&s_n&\end{array} \right)\] 

Adding up the rows and columns shows $\beta_r$ projects to $\sum \limits_{k = 1}^n x_k s_k$ and $\sum \limits_{k=1}^n y_k t_k$. Define 
\[\beta = \sum \limits_{r \in R(U)}\beta_r\]

Only finitely many of the terms of this sum are nonzero, and $\beta$ indeed maps to $\alpha$. 

\end{proof}

\begin{remark}\label{counterexamplewisedef}

The map $A^{S \times_R T} \rightarrow A^S \times_{A^R} A^T$ is \textit{not} an isomorphism, in general. A counterexample is found already when $E = (Sets)$, $A = \Q$ and $R = \{r\}$. 

Consider $S := \{x, y\}$ and $T := \{x', y'\}$ with their unique maps to $R$. Then $\Q^{S \times_R T} \rightarrow \Q^S \times_{\Q^R} \Q^T$ is surjective but not injective. For example, $(x, y) - (x, y') + (x', y') - (x', y)$ goes to zero. Hence the functor $S \mapsto A^S$ needn't commute with finite limits and is not left exact. 

Applying $\text{Sym}$ to the above counterexample shows the free algebra functor $S \mapsto A[S]$ isn't left exact either, contradicting a claim made in \cite{wisedef1}. Forthcoming work will show the conclusions in \cite{wisedef1} which rest on this erroneous claim remain true. For their proof, an analogue of Lemma \ref{finlimfreecovers} suffices. 

\end{remark}

\begin{remark}\label{freepullbackisasheaf}
Because the functor $\lambda : E \rightarrow \Amod$ of Lemma \ref{finlimfreecovers} is not left exact, it doesn't induce a morphism of sites in the other direction. It is cocontinuous nonetheless, inducing a morphism of sites
\[E \rightarrow \Amod\]
The left exact left adjoint belonging to this morphism is precisely $F \mapsto (\lambda^*F)^{sh}$, the sheafification of precomposition by $\lambda$. We note that $\lambda^*F$ is already a sheaf. 

For $\{S_i \rightarrow T\}$ a cover in $E$, the projections $A^{S_i \times_T S_j} \rightarrow A^{S_i}$ factor through the product $A^{S_i} \times_{A^T} A^{S_j}$. Taking sections over $F$, we get a sequence
\[F(A^T) \rightarrow \prod F(A^{S_i}) \rightrightarrows \prod F(A^{S_i} \times_{A^T} A^{S_j}) \hookrightarrow \prod F(A^{S_i \times_T S_j})\]

The last map is injective because the map $A^{S_i \times_T S_j} \rightarrow A^{S_i} \times_{A^T} A^{S_j}$ is covering. The sheaf condition for $\lambda^*F$ is that the diagram formed by the first arrow and the composites of the pair of arrows with the injection should be an equalizer. However, the sheaf condition on $\Amod$ ensures that the diagram without the final injection is an equalizer. Postcomposing by an injection preserves such an equalizer diagram. 

\end{remark}

\begin{proposition}\label{rhoisequiv}
The functor $\rho : \EExt_A^1(\inpt, K) \rightarrow B \hh{K}$ of Defintion \ref{definerho} is an equivalence. 
\end{proposition}

\begin{proof}

In the process of defining $\rho$, we remarked that it's fully faithful. It remains to show essential surjectivity. 

For any $N \in \Amod$, write $j_N : \Amod/N \rightarrow \Amod$ for the localization morphism of topoi. Write $\lambda : E \rightarrow \Amod$ for the functor $S \mapsto A^S$ and $\lambda^*$ for the induced functor on sheaves $F \mapsto F \circ \lambda$. Write $\alpha_N$ for the map $A \times N^{\times 2} \rightarrow N$ which is $(a, n, n') \mapsto a.(n+n')$ on sections. 

Let $P$ be an $\hh{K}|_M = \hh{K \oplus M}$-torsor on $\Amod/M$. Let $\{M_i \rightarrow M\}_I$ be a cover on which $P$ is trivial. For $N \in \Amod/M$, define

\[L_N := \lambda^* j_{N!} P|_N\]

This is the sheafification of the functor $U \mapsto \bigsqcup \limits_{f \in \Hom_A(A^U, N)} P|_N(f)$. 

Write $L_i := L_{M_i}$ and $L := L_M$ for brevity. There are maps $p_N : L_N \rightarrow N$ sending $P(f)$ to the section of $N$ corresponding to $f$. If $N \rightarrow N' \in \Amod/M$, then $L_N = L_{N'} \times_{N'} N$ and $p_N$ is the projection. 

We will show $p_M : L \rightarrow M$ fits into an extension of modules which maps to $P$ under $\rho$. First, we must augment the sheaf of sets $L$ with an $A$-module structure.  

Define $R \subseteq \Amod/M$ to be the full subcategory on those $N$ whose structural morphism to $M$ factors through some $M_i$ (the sieve generated by the cover). We will produce an $A$-module structure on $L_N$ for each $N \in R$ and then descend to $L$. 

Remark that $\lambda^* j_{N!} \hh{N' | N} = N'$ functorially in $N'$. Choose an $\hh{K}$-isomorphism $f : P|_N \simeq \hh{K}|_N = \hh{K \oplus N}$. Give $L_N$ the induced $A$-module structure from the $K$-isomorphism $\lambda^* j_{N!} f : L_N \simeq K \oplus N$. We claim this $A$-module structure is independent of the choice of $f$. 

To that end, let $g : P|_N \simeq \hh{K \oplus N}$ be another $\hh{K}$-isomorphism. For the map $g \circ f^{-1}$ to be an $\hh{K}$-equivariant map of representable sheaves on $\Amod/N$, it must come from a map of $A$-modules. Since $g \circ f^{-1}$ comes from an $A$-module homomorphism, $f$ and $g$ endow $L_N(U)$ with the same $A$-module structure. 

Since the $A$-module structure on each is well-defined, the equality of sheaves $L_N = L_{N'} \times_{N'} N$ is promoted to one of modules. In particular, the projection maps $L_N \rightarrow L_{N'}$ are each $A$-module maps.

We want to construct $\alpha_M$ using the cover. By the definition of the topology on $\Amod$, if $\{M_i \rightarrow M\}_I$ is a cover, then $\{M_i^{\times 2} \rightarrow M^{\times 2}\}_I$ is also a cover. It follows that the pullback $\{A \times L_i^{\times 2} \rightarrow A \times L^{\times 2}\}$ is covering in $E$. 

Since the topology on $E$ is subcanonical, 
\[\hh{L}(A \times L^{\times 2}) \rightarrow \prod \hh{L}(A \times L_i^{\times 2}) \rightrightarrows \prod \hh{L}(A \times (L_i \times_L L_j)^{\times 2})\]
is an equalizer. The commutativity of the following diagram
\begin{equation}\label{modulemapsheaf}
\begin{tikzcd}
A \times (L_i \times_L L_j)^{\times 2} \arrow{r} \arrow{d}{\alpha_{L_i \times_L L_j}} \arrow[phantom, "\circ"]{dr}       &A \times L_i^{\times 2} \arrow{r} \arrow{d}{\alpha_{L_i}}         &A \times L^{\times 2} \arrow[dashed]{d}       \\
L_i \times_L L_j \arrow{r}       &L_i \arrow{r}        &L
\end{tikzcd}
\end{equation}

\noindent ensures the existence of the dashed arrow. Commutativity is the statement that the projections $L_i \times_L L_j = L_{M_i \times_M M_j} \rightarrow L_i$ are $A$-module maps. 

The dashed arrow $A \times L^{\times 2} \dashrightarrow L$ \textit{defines} addition and scalar multiplication for $L$. Since $\hh{L}$ is a sheaf, equality between two arrows to $L$ may be checked after pulling back along a cover of $M$. This guarantees commutativity, associativity, etc.

The epimorphism $\bigoplus L_i \rightarrow \bigoplus M_i \rightarrow M$ factors through $L$, guaranteeing $L \rightarrow M$ to be epimorphic. The kernel is seen to be $K$ by pulling back $L \rightarrow M$ along any $M_i \rightarrow M$. 

It remains to show the extension 
\[\xi : 0 \rightarrow K \rightarrow L \rightarrow M \rightarrow 0\]
represents $P$; that $\rho(\xi) \simeq P$. We build an isomorphism for any $N$ in the sieve $R$ and show it's independent of choices. Choose two $\hh{K}$-isomorphisms $f, g : P|_N \simeq \hh{K}|_N$. Apply $\hh{(\lambda^* j_{N!} \inpt)|N}$ to both and form the commutative diagram: 
\[\begin{tikzcd}[column sep=small, row sep=small]
        &\hh{K}|_N \arrow[no head, bend left, "\sim f"]{dr}      &       \\
\hh{L_N|N} \arrow[no head, bend left, "\sim f"]{ur} \arrow[no head, bend right, swap, "\sim g"]{dr} \arrow[dashed]{rr}      &       &P|_N       \\
        &\hh{K}|_N \arrow[no head, bend right, swap, "\sim g"]{ur}
\end{tikzcd}\]

This verifies compatibility of the locally defined isomorphisms $\hh{L_N|N} \simeq P|_N$ and we obtain a global $\hh{L|M} \simeq P$. Hence $\rho$ is essentially surjective.

\end{proof}

\begin{remark}\label{extdefstack}
Lemma \ref{rhoisequiv} implies that $\EExt_A^1(\inpt, K) \rightarrow \Amod$ is a stack, in fact a form of the trivial $\hh{K}$-gerbe. The proof checked descent by relying heavily on $\rho$. 

In the same way, $\EExt_{A'}^1(\inpt, K) \rightarrow \AAmod$ is an $\hh{K}$-gerbe. However, $\pi_*\EExt_{A'}^1(\inpt, K) \rightarrow \Amod$ is a stack but no longer an $\hh{K}$-gerbe.

Diagram \eqref{definingdef} defining $\DDef(\inpt, u, K)$ describes it as a strict fiber product of stacks over $\Amod/M$, so it's also a stack. 
\end{remark}

\begin{remark}
The topology was used in the proof only to the extent that, 
if $\{M_i \rightarrow M\}_I$ is a cover, then $\{M_i^{\times 2} \rightarrow M^{\times 2}\}_I$ is also a cover in $E$. We could vary the topology so that the $\Lambda$ in the definition of the topology could only have at most two elements, and the proof would still work. 

We speculate that allowing $\Lambda$ to have at most three elements would suffice for $\EExt^2_A(\inpt, K)$ to form a 2-gerbe, and consider this infinite hierarchy of topologies curious. 
\end{remark}

We finish the section with a few more basic properties of the site $\Amod$. Define $\mathcal{P}$ as the presheaf on $\Amod$ (resp. define a sheaf $\mathcal{P}_E$ on $E$) whose value on $M$ is the set of submodules (resp. subsheaves) of $M$. Precomposing by the forgetful functor $\Amod \rightarrow E$, we regard $\mathcal{P}_E$ as a sheaf on $\Amod$ and $\mathcal{P}$ as a subpresheaf. 

Since $E$ has a set of generators, $\mathcal{P}_E(M)$ and $\mathcal{P}(M)$ are indeed sets. Restriction maps are pullbacks of subobjects. 

\begin{lemma}
The presheaf $\mathcal{P}$ on $\Amod$ is a sheaf. 
\end{lemma}

\begin{proof}

Let $\{M_i \rightarrow M\}_I$ be a cover, with submodules $N_i \subseteq M_i$. Write $M_{ij} := M_i \times_M M_j$. Suppose the pullbacks $N_i|_{M_{ij}} = N_j|_{M_{ij}} \subseteq M_{ij}$ are equal. We want to exhibit a submodule $N \subseteq M$ whose pullbacks to each $M_i$ are precisely $N_i$. 

Since $\mathcal{P}_E$ is a sheaf, the above descent data furnishes a subsheaf of sets $N \subseteq M$ on $E$; we must endow $N$ with a submodule structure. We get a diagram as in \eqref{modulemapsheaf} by replacing $L$ by $N$, and the same argument produces the submodule structure. 

\end{proof}

\begin{corollary}\label{arrstack} The arrow category $q : \Arr(\Amod) \rightarrow \Amod$ \cite[3.15]{vistolistacks} is a stack, the functor sending an arrow to its codomain.
\end{corollary}

\begin{proof}

Isomorphisms form a sheaf for $\Arr(\Amod)$ because $\IIsom_A(N, N')$ is a sheaf. 

Any arrow $N \rightarrow M \in \Amod$ factors as $N \twoheadrightarrow P \hookrightarrow M$, an epimorphism composed with a monomorphism. Considering $\mathcal{P}$ as a fibered category, factor the functor $q$ as 
\[\Arr(\Amod) \rightarrow \mathcal{P} \rightarrow \Amod\]
The first arrow sends $N \rightarrow M$ to the image $P \subseteq M$. Since $\mathcal{P}$ is a sheaf, we need only show $\Arr(\Amod) \rightarrow \mathcal{P}$ satisfies descent for the induced topology \cite[06NU, 09WX]{sta}. The corollary follows. 

The induced topology refers to cartesian arrows over a cover in the base site. In other words, a cover in $\mathcal{P}$ is a cover $\{M_i \rightarrow M\}_I$ in $\Amod$ together with a choice of subobject $N \subseteq M$ and its pullbacks to $M_i$. 

Descent data for $\Arr(\Amod)$ here refers to a choice of epimorphism $M'_i \twoheadrightarrow N|_{M_i}$, isomorphisms between the pullbacks of $M'_i$ and $M'_j$ along $M_{ij}$'s two projections compatible with the epimorphisms, and compatibility of those isomorphisms on $M_{ijk}$. Remark that the kernel of each epimorphism must be the same, say $K$. This is precisely a descent datum for $\EExt_A^1(\inpt, K) \rightarrow \Amod$, necessarily effective by Remark \ref{extdefstack}. We obtain an epimorphism $M' \twoheadrightarrow N$, also with kernel $K$, which pulls back to each $M'_i \twoheadrightarrow N|_{M_i}$ and verifies descent.

\end{proof}

We have now developed enough technology to solve the deformation problem.

\section{Cohomology on A-mod}\label{deformations}

We can quickly solve the deformation problem with an algebraic statement. This theorem yields an obstruction in degree-two cohomology. The rest of the section is devoted to the proof of Theorem \ref{newextiscohom}.

\openrefcount{defisagerbe}
\begin{theorem}\label{newdefisagerbe}
The fibered category $\DDef(\inpt, u, K) \rightarrow \Amod$ is a gerbe banded by $\hh{K}$. 
\end{theorem}
\closerefcount

\begin{proof}

We've seen already in Remark \ref{extdefstack} that $\DDef(\inpt, u, K)$ is a stack. 

Given an automorphism 
\[\begin{tikzcd}
0 \arrow{r}       &K \arrow{r} \arrow[equals]{d}      &M' \arrow{r} \arrow{d}{\sim}     &M \arrow{r} \arrow[equals]{d}      &0      \\
0 \arrow{r}       &K \arrow{r}      &M' \arrow{r}     &M \arrow{r}      &0
\end{tikzcd}\]
of a global section of $\DDef(\inpt, u, K) \rightarrow \Amod/M$, subtract the identity. The resulting morphism of chain complexes is zero on $K$ and $M$, and the map $M' \rightarrow M'$ factors through $M$ and $K$. Automorphisms of deformations are thereby in bijection with $\hh{K}(M)$.

It remains to show $\DDef(\inpt, u, K)$ is locally nonempty and two sections are locally isomorphic. 

For both, we may assume $M=A^S$ for some sheaf of sets $S \in E$ by localizing in $\Amod$. Tensor the short exact sequence \eqref{fixalgdef} by $\otimes A'^S$ to get a canonical deformation of $id_{J \otimes A^S}$: 
\[\underline{\alpha} : \quad 0 \rightarrow J \otimes A^S \rightarrow A'^S \rightarrow A^S \rightarrow 0\]

To deform an arbitrary map $u : J \otimes M \rightarrow K$, simply localize in $M$ and pushout $\underline{\alpha}$ by $u$. Observe $\theta(u \smile \underline{\alpha}) = u$.

Now we show sections of $\DDef(\inpt, u, K)$ are locally isomorphic. Choose an extension 
\[\xi : \quad 0 \rightarrow K \rightarrow M' \rightarrow M \rightarrow 0\]
with $\theta(\xi) = u$. Since extensions are locally trivial, we may choose a cover $\{A^S \rightarrow M\}$ so that each $S$ lifts to $M'$. We obtain a morphism of extensions 

\[\begin{tikzcd}
0 \arrow{r}      &J \otimes A^S \arrow{r} \arrow[dashed, "u"]{d}      &A'^S \arrow{r} \arrow[dashed]{d}      &A^S \arrow{r} \arrow{d}       &0      \\
0 \arrow{r}      &K \arrow{r}      &M'\arrow{r}     &M \arrow{r}     &0
\end{tikzcd}\]
witnessing that $\xi|_{A^S} \simeq u \smile \underline{\alpha}$. The induced map on the kernel is forced to be $u$.

\end{proof}

Now that we have a degree-two cohomological obstruction, we must work explicitly with the cohomology groups $H^p(\Amod/M, \hh{K})$. The remainder of the section proves Theorem \ref{extiscohom}.

\begin{lemma}\label{injdetectexact}
A complex of $A$-modules 
\[C_\bullet : \quad     \cdots \rightarrow C_{p+1} \rightarrow C_p \rightarrow C_{p-1} \rightarrow \cdots\]
is exact if and only if, for any injective $A$-module $K$, the complex of homomorphisms into $K$ is: 
\[\hh{K}(C_\bullet) : \quad     \cdots \leftarrow \Hom_A(C_{p+1}, K) \leftarrow \Hom_A(C_p, K) \leftarrow \Hom_A(C_{p-1}, K) \leftarrow \cdots\]
\end{lemma}

\begin{proof}

The ``only if'' is clear. For the reverse implication, we want to show the map $C_{p+1} \rightarrow \ker d_p$ is surjective. It suffices to show the cokernel $(\ker d_p)/C_{p+1} = 0$. Embed the quotient into an injective $A$-module $K$: $(\ker d_p)/C_{p+1} \subseteq K$. This corresponds to $\ker d_p \rightarrow K$ such that $C_{p+1} \rightarrow \ker d_p \rightarrow K$ is zero. By injectivity of $K$, we get a factorization $\ker d_p \rightarrow C_p \dashrightarrow K$ such that $C_{p+1} \rightarrow C_p \rightarrow K$ is zero. 

By the hypothesized exactness of $\hh{K}(C_\bullet)$, we get a factorization $C_p \rightarrow C_{p-1} \dashrightarrow K$. Then our original map factors $\ker d_p \rightarrow C_p \rightarrow C_{p-1} \rightarrow K$, and the first composition is zero. Thus the induced map $(\ker d_p)/C_{p+1} \rightarrow K$ is zero, but it was assumed to be an embedding. Therefore $(\ker d_p)/C_{p+1} = 0$, the map $C_{p+1} \rightarrow \ker d_p$ is surjective, and the complex has no cohomology in the $p$th degree.

\end{proof}

\begin{proposition}\label{injvanishing}
Given an injective $A$-module $K$, the higher cohomology of $K$ all vanishes. That is, $H^p(\Amod/M,\hh{K}) = 0$ for $p \geq 1$.
\end{proposition}

Abstract properties of derived functors turn our main theorem into an immediate consequence of the previous proposition. We show how before providing the proof, the most complicated in this paper. 

\openrefcount{extiscohom}
\begin{theorem}\label{newextiscohom}
The groups of $p$-extensions are all equivalent to cohomology of $\hh{K}$ on the site $\Amod$ : $\Ext^p_A(M, K) \simeq H^p(\Amod/M, \hh{K})$. 
\end{theorem}
\closerefcount

\begin{proof}

Proposition \ref{injvanishing} shows that $H^p(\Amod/M, \hh{K})$ is a universal $\delta$-functor in $K$. Since $H^0(\Amod/M, \hh{K}) := \Hom_A(M, K)$ and $\Ext^p_A(M, K)$ is defined to be a universal $\delta$-functor in $K$, we get a unique isomorphism $H^p(\Amod/M, \hh{K}) \simeq \Ext^p_A(M, K)$ of $\delta$-functors by \cite[III.1.2.1]{hartshorne}. 

\end{proof}

\begin{proof}[Proof of Proposition \ref{injvanishing}]

We will prove exactness of the \v{C}ech Complex in a series of lemmas to follow. We recall a well-known reduction to the vanishing of \v{C}ech Cohomology in the meantime (\cite[01EV]{sta}, usually attributed to Cartan), as we will need the details in Lemma \ref{finitegencech}. 

Assume inductively that $H^i(M, \hh{K}) = 0$ for $0 < i < p$ and any injective $K$. Proposition \ref{rhoisequiv} yields the base case: 
\[H^1(M, \hh{K}) = \Ext_A^1(M, K) = 0\]
for $K$ injective. Consider a cover $\{M_i \rightarrow M\}_I$ and an injective $A$-module $K$. 

The \v{C}ech Spectral Sequence \cite[V.3.3]{sga4} is: 
\[H^j(M_\bullet, \underline{H}^k\hh{K}) \Rightarrow H^{j+k}(M, \hh{K})\]
Here $\underline{H}^k \hh{K}$ is the presheaf $N \rightarrow H^k(N, \hh{K})$ and is zero for $0 < k < p$ by inductive assumption. The only possibly nonzero terms on the diagonal $j + k = p$ are $H^p(M_\bullet, \underline{H}^0 \hh{K})$ and $H^0(M_\bullet, \hh{K})$. The filtration on degree $p$ cohomology is expressed by the exact sequence 
\begin{equation}\label{filtrationoncech}
0 \rightarrow H^p(M_\bullet, \underline{H}^0 \hh{K}) \rightarrow H^p(M, \hh{K}) \rightarrow H^0(M_\bullet, \underline{H}^p \hh{K}) \rightarrow \cdots
\end{equation}
The \v{C}ech spectral sequence and therefore this short exact sequence are natural with respect to refinement of the cover. The map on the right arises from the restriction of the presheaf $\underline{H}^p$. 

In order to show $H^p(M, \hh{K})$ vanishes, pick an element $\alpha$. Then \cite[01FW]{sta} allows us to choose a cover $\{M_i \rightarrow M\}_I$ so that $\alpha|_{M_i} = 0$, so $\alpha \in H^p(M_\bullet, \underline{H}^0 \hh{K})$ by the exactness of \eqref{filtrationoncech}. It suffices therefore to show \v{C}ech Cohomology vanishes. 

To that end, fix a total ordering on $I$. Write $M_{i_0 \cdots i_p} := M_{i_0} \times_M \cdots \times_M M_{i_p}$. We get a semi-simplicial object whose $p$th simplices are $\bigoplus M_{i_0 \cdots i_p}$, the sum ranging over ordered $(p+1)$-tuples of indices $i_0 \leq i_1 \leq \cdots \leq i_p$. The $j$th face map projects away from $i_j$. Take alternating sums to obtain a complex: 
\begin{equation}\label{cechcomplexmodules}
\cdots \rightarrow \bigoplus M_{i_0 \cdots i_{p+1}} \rightarrow \bigoplus M_{i_0 \cdots i_p} \rightarrow \bigoplus M_{i_0 \cdots i_{p-1}} \rightarrow \cdots \rightarrow \bigoplus M_i \rightarrow M \rightarrow 0
\end{equation}
\v{C}ech Cohomology results from applying $\hh{K}$ to this sequence and taking cohomology. Complex \eqref{cechcomplexmodules} is exact precisely when \v{C}ech Cohomology vanishes by Lemma \ref{injdetectexact}. We've reduced the proof to the following Lemma \ref{completetheproof}. 

\end{proof}

\begin{lemma}\label{completetheproof}
Complex \eqref{cechcomplexmodules} is exact. 
\end{lemma}

We prove this lemma after first handling a few special cases. 

\begin{lemma}\label{singleeltcech}
Case 1 of Lemma \ref{completetheproof}: If $\{T \rightarrow M\}$ is a cover consisting of a single element (i.e. $I = \{*\}$), then the complex \eqref{cechcomplexmodules} is exact. 
\end{lemma}

\begin{proof}

Remark in particular that $T \rightarrow M$ is an epimorphism. Write $S$ for its kernel. We describe maps of sheaves of modules on sections $t \in \Gamma(U, T)$ and $s_i \in \Gamma(U, S)$ over some $U \in E$ for convenience. 

The shearing map $T \oplus S^{\oplus p} \overset{\sim}{\rightarrow} T \times_M \cdots \times_M T$ sending $(t, s_1, \cdots s_p)$ to the partial sums $(t, t + s_1, t + s_1 + s_2, \cdots, t + \Sigma s_i)$ is an isomorphism. Under this isomorphism, the semi-simplicial module yielding \eqref{cechcomplexmodules} has $p$-simplices $T \oplus S^{\oplus p}$ and face maps given by

\[d_i(t, s_1, \cdots, s_p) := \left\{ 
\begin{tikzcd}
(t, s_1, \cdots, s_{i-1}, s_i + s_{i+1}, s_{i+2}, \cdots, s_p)     &\text{if }i \neq 0, p     \\
(t + s_1, s_2, \cdots, s_p)       &i=0        \\
(t, s_1, \cdots, s_{p-1})         &i=p
\end{tikzcd}
\right.\]

The reader is invited to verify the semi-simplicial axiom and verify this assignment yields an isomorphism of semi-simplicial modules with $T \times_M \cdots \times_M T$. Witness the similarity to the simplicial construction of $EG$ in \cite[pg. 128]{concisealgtop}. 

Now we check by hand that the normalized chain complex associated to $T \oplus S^{\oplus p}$ is exact. The normalized chain complex in degree $p$ is the intersection of all the kernels of the $d_i$, for $i \neq 0$ -- its differentials are precisely $d_0$. Consider a local section $(t, s_1, \cdots, s_p)$ of the $p$-th degree of the normalized chain complex. 

In order to be in the kernel of $d_i$ for $i \neq 0, p$, $t = 0$, $s_i = -s_{i+1}$ and all $s_j = 0$ except $j = i, i+1$. Consider a few cases: 
\begin{itemize}
\item $p \geq 3$: Varying $i$ implies $t = s_1 = \cdots = s_p = 0$. 
\item $p = 2$: For $d_2$ to vanish we must also have $s_1 = 0$, and again $t = s_1 = s_2 = 0$. 
\item $p = 1$: For $d_1$ to vanish, $t = 0$. 
\item There are no requirements for $p = 0$. 
\end{itemize}

The augmented normalized chain complex is thereby seen to be 
\[0 \rightarrow S \rightarrow T \rightarrow M \rightarrow 0\]
with the natural maps. This is exact by assumption. The normalized chain complex is well known \cite[I.1.3.3]{illusie1} to be quasi-isomorphic to the unnormalized chain complex \eqref{cechcomplexmodules}.

\end{proof}

\begin{remark}
We caution the reader that $\bigoplus$ and products over $M$ do not commute, and hence \eqref{cechcomplexmodules} is not a series of fiber products of $\bigoplus M_i \rightarrow M$. That is, the problem does not reduce entirely to Lemma \ref{singleeltcech}. 

\end{remark}

\begin{lemma}\label{finitegencech}
Case 2 of Lemma \ref{completetheproof}: Suppose $M$ has a finite set of global sections $\Lambda \subseteq \Gamma(E, M)$ so that the induced map from the free module on the constant sheaf $A^\Lambda \rightarrow M$ is covering. Then complex \eqref{cechcomplexmodules} is exact. 
\end{lemma}

\begin{proof}

The covering condition allows us to localize in $E$ so that $\Lambda$ lifts to some $M_i$, say $M_0$. Then the hypothesized cover factors as $A^\Lambda \rightarrow M_0 \rightarrow M$ and $M_0 \rightarrow M$ is a cover.

The inclusion $M_0 \subseteq \{M_i\}_I$ is a refinement of covers. The short exact sequence \eqref{filtrationoncech} is contravariant under refinements: 
$$\begin{tikzcd}
0 \arrow{r}         &H^p(M_\bullet, \underline{H}^0\hh{K}) \arrow{d} \arrow{r}       &H^p(M, \hh{K}) \arrow[equals]{d} \arrow{r}     &\cdots      \\
0 \arrow{r}         &H^p(\{M_0\}, \underline{H}^0\hh{K}) \arrow{r}       &H^p(M, \hh{K}) \arrow{r}      &\cdots
\end{tikzcd}$$
By Lemma \ref{singleeltcech}, the group $H^p(\{M_0\}, \underline{H}^0\hh{K})$ vanishes, but the injection $H^p(M_\bullet, \underline{H}^0\hh{K}) \hookrightarrow H^p(M, \hh{K})$ factors through this group; this implies $H^p(M_\bullet, \underline{H}^0\hh{K}) = 0$. Equivalently, $\hh{K}(M_\bullet)$ is exact. Since $K$ was any injective module, $\eqref{cechcomplexmodules}$ is exact by Lemma \ref{injdetectexact}.

\end{proof}

We are finally ready to complete the proof of Lemma \ref{completetheproof}

\begin{proof}[Proof of Lemma \ref{completetheproof}]

We often use the observation that, in order to verify exactness of the sequence of sheaves \eqref{cechcomplexmodules}, we may freely localize in $E$. 

To show \v{C}ech Cohomology vanishes, suppose some section $\beta = \sum m_j \in \Gamma(U, \bigoplus M_{i_0 \cdots i_p})$ maps to zero. Localize so that $\beta$ is a global section. Define $N$ as the image of the map $A^{\{m_j\}} \rightarrow M$ adjoint to the map from the constant sheaf $\{m_j\} \subseteq M_{i_0 \cdots i_p} \rightarrow M$. 

Write $N_i := M_i \times_M N$ and form the following diagram: 
$$\begin{tikzcd}
\cdots \arrow{r}      &\bigoplus N_{i_0 \cdots i_{p+1}} \arrow[hook]{d} \arrow{r} \arrow[very near start, phantom, "\ulcorner"]{dr}      &\bigoplus N_{i_0 \cdots i_p} \arrow[hook]{d} \arrow{r} \arrow[very near start, phantom, "\ulcorner"]{dr}    &\bigoplus N_{i_0 \cdots i_{p-1}} \arrow[hook]{d} \arrow{r} \arrow[very near start, phantom, "\ulcorner"]{dr}      &\cdots \arrow{r}     &N \arrow[hook]{d}     \\
\cdots \arrow{r}      &\bigoplus M_{i_0 \cdots i_{p+1}} \arrow{r}      &\bigoplus M_{i_0 \cdots i_p} \arrow{r}     &\bigoplus M_{i_0 \cdots i_{p-1}} \arrow{r}       &\cdots \arrow{r}     &M      
\end{tikzcd}$$
Each map $N_{i_0 \cdots i_q} \rightarrow M_{i_0 \cdots i_q}$ is the pullback of $N \hookrightarrow M$, so the vertical arrows are monomorphisms. 

By construction, there is a preimage $\tilde{\beta} \in \Gamma(E, \bigoplus N_{i_0 \cdots i_p})$ of $\beta$. Moreover, $\tilde{\beta}$ maps to zero under the differential by injectivity of the vertical maps. The module $N$ was defined as the image of $A^{\{m_j\}}$, so it falls under the jurisdiction of Lemma \ref{finitegencech}, and the top row is exact. Then $\tilde{\beta}$ is a boundary. This concludes the proof. 

\end{proof}

Armed with the isomorphism of Theorem \ref{extiscohom}, we now undertake its study.

%\section{Extensions and Cohomology}

%\subfile{2cechextandcohom}

\section{Extensions and Cohomology}\label{extandcohom}

%Ext and Cohomology

This section describes the isomorphisms of Theorem \ref{extiscohom} in degrees $p = 1, 2$. We use this description to prove Theorem \ref{obstcomparison}. 

\begin{lemma}\label{explicitp1}
The isomorphism $\Ext_A^p(M, K) \simeq H^p(\Amod/M, \hh{K})$ of Theorem \ref{extiscohom} in degree $p = 1$ is the restriction of the functor $\rho$ of Definition \ref{definerho} to isomorphism classes. 
\end{lemma}

\begin{proof}

Given a short exact sequence 

\[\gamma : 0 \rightarrow K \rightarrow N' \rightarrow N \rightarrow 0,\]

The diagram with horizontal arrows the boundary maps for the long exact sequences of $\Ext_A^p$ and $H^p$ 

\[\begin{tikzcd}
\Hom_A(M, N) \arrow{r} \arrow[equals]{d}        &\Ext_A^1(M, K) \arrow{d}{\rho}     \\
H^0(\Amod/M, \hh{N}) \arrow{r}        &H^1(\Amod/M, \hh{K})
\end{tikzcd}\]

\noindent commutes. Indeed, a map $M \rightarrow N$ is sent via the boundary map for $H^p$ to the $\hh{K}$-torsor of sections of $N' \rightarrow N$ (see \ref{torsorgerbesignconvention}): 
\[P(T) := \left\{ 
\begin{tikzcd}
        &       &N' \arrow{d}     \\
T \arrow[dashed]{urr} \arrow{r}     &M  \arrow{r}       &N
\end{tikzcd}
\right\}\]

The boundary map for $\Ext_A^p$ sends the map $M \rightarrow N$ to the pullback $\gamma|_M$ of the extension along the map. Under $\rho$, this extension is sent to the $\hh{K}$-torsor represented by $M \times_N N'$ in $\Amod/M$. Sections of $M \times_N N' \rightarrow M$ and elements of $P$ are in a canonical bijection which respects the action of $\hh{K}$.

Suppose now that $N'$ is an injective $A$-module. Then the horizontal boundary maps are epimorphisms, and we see that $\rho$ is the same as the isomorphism provided by Theorem \ref{extiscohom}. 

\end{proof}

We must now develop a considerable amount of technology to deal with the $p=2$ case. Fix notation for two 2-extensions for the rest of the section: 

\begin{equation}\label{2extsetup}
\begin{tikzcd}
\xi:      &0 \arrow{r}      &K \arrow{r}       &X \arrow{r}      &Y \arrow{r}      &M \arrow{r}      &0      \\
\eta:     &0 \arrow{r}      &K \arrow{r}      &X' \arrow{r}     &Y' \arrow{r}     &M \arrow{r}      &0
\end{tikzcd}
\end{equation}

Write $P$ for the module $\coker(K \rightarrow X) \simeq \ker(Y \rightarrow M)$ and $P'$ likewise for $\coker(K \rightarrow X') \simeq \ker(Y' \rightarrow M)$. 

Define the trivial 2-extension as 
\[\begin{tikzcd}
\underline{0} :         &0 \arrow{r}      &K \arrow[equals]{r}      &K \arrow{r}{0}      &M \arrow[equals]{r}     &M \arrow{r}      &0
\end{tikzcd}\]

\begin{definition}
A \textit{butterfly} $\xi \simeq \eta$ between two 2-extensions is a completion of \eqref{standardbutterfly} or of the equivalent diagram \eqref{diamondbutterfly}. They form the isomorphisms in a 2-groupoid $\EExt_A^2(M, K)$. The 2-isomorphisms are given by isomorphisms of the completions $Q \simeq Q'$ which commute with all of the structure maps. 

We owe the concept to \cite{butterflies} and \cite{sga7i}. As in the latter, we'll be concerned only with the abelian case. Some background on butterflies in our context is recalled in Appendix \ref{butterfliesappendix}. 

The fibered category $\IIsom(\xi, \eta) \rightarrow \Amod/M$ has fiber over $N \rightarrow M$ given by the category of butterflies between $\xi|_N$ and $\eta|_N$. Write $\Split(\xi)$ for $\IIsom(\xi, \underline{0})$. 
\end{definition}

\begin{equation}\label{standardbutterfly}
\begin{tikzcd}
        &0 \arrow{dr}      &       &       &&0      &       \\
0 \arrow{r}      &K \arrow{r} \arrow[equals]{dd}     &X\arrow[dashed]{dr} \arrow{rr}     &      &Y \arrow{r} \arrow{ur}       &M  \arrow[equals]{dd}\arrow{r}    &0      \\
        &       &       &Q \arrow[dashed]{dr} \arrow[dashed]{ur}      &       &       &       \\
0 \arrow{r}      &K \arrow{r}     &X'\arrow{rr} \arrow[dashed]{ur}     &      &Y' \arrow{r} \arrow{dr}      &M  \arrow{r}    &0       \\
        &0 \arrow{ur}      &       &       &&0      &
\end{tikzcd}
\end{equation}

\begin{equation}\label{diamondbutterfly}
\begin{tikzcd}
&   &0 \arrow{dr}   &   &0   &   \\
&0 \arrow{dr}        &       &P \arrow{dr} \arrow{ur}     &       &0       \\
0 \arrow{dr}    &        &X \arrow{ur} \arrow[dashed]{dr}     &       &Y \arrow{dr} \arrow{ur}     &       &0     \\
&K \arrow{ur} \arrow{dr}       &       &Q \arrow[dashed]{ur} \arrow[dashed]{dr}      &       &M \arrow{ur} \arrow{dr}      \\
0 \arrow{ur}     &        &X' \arrow[dashed]{ur} \arrow{dr}     &       &Y' \arrow{ur} \arrow{dr}     &       &0      \\
&0 \arrow{ur}        &       &P' \arrow{ur} \arrow{dr}     &       &0     \\
&   &0 \arrow{ur}   &   &0   &
\end{tikzcd}
\end{equation}

%%%butterfly without all the zeros
\begin{comment}
\begin{equation}
\begin{tikzcd}
        &       &P \arrow{dr}     &       &       \\
        &X \arrow{ur} \arrow[dashed]{dr}     &       &Y \arrow{dr}     &       \\
K \arrow{ur} \arrow{dr}       &       &Q \arrow[dashed]{ur} \arrow[dashed]{dr}      &       &M      \\
        &X' \arrow[dashed]{ur} \arrow{dr}     &       &Y' \arrow{ur}     &       \\
        &       &P' \arrow{ur}     &       &       
\end{tikzcd}
\end{equation}
\end{comment}

In the first diagram, the NW-SE and SW-NE diagonals in the interior are short exact sequences. In the diamond-shaped diagram, each line is a short exact sequence. We show that $\IIsom(\xi, \eta)$ is an $\hh{K}$-banded gerbe even though we are particularly interested in $\Split(\xi)$.

\begin{lemma}
The fibered category $\IIsom(\xi, \eta) \rightarrow \Amod/M$ is a stack. 
\end{lemma}

\begin{proof}

All the data are local, so $\IIsom(\xi, \eta)$ must be a stack. To show descent data are effective, one can locally

\begin{itemize}

\item Build an arrow $Q \rightarrow M$ ($\Arr(\Amod/M)$ is a stack -- Corollary \ref{arrstack}). 

\item Build factorizations $Q \rightarrow Y \rightarrow M$ and $Q \rightarrow Y' \rightarrow M$ ($\hh{M}$ is a sheaf).

\item Check exactness of the short exact sequences 
\[0 \rightarrow X' \rightarrow Q \rightarrow Y \rightarrow 0\]
and 
\[0 \rightarrow X \rightarrow Q \rightarrow Y' \rightarrow 0\]
(the composite $\EExt_A(\inpt, X') \rightarrow \Amod/Y \overset{j_!}{\rightarrow} \Amod/M$ is a stack, if $j$ is the localization morphism of topoi). 

\item Check commutativity of the North, West, and South diamonds in Diagram \eqref{diamondbutterfly} ($\HHom_A(X, Y)$ is a sheaf on $\Amod/M$). 

\end{itemize}

\end{proof}

\begin{lemma}\label{isomisagerbe}
The stack $\IIsom(\xi, \eta) \rightarrow \Amod/M$ is a gerbe banded by $\hh{K}$. 
\end{lemma}

\begin{proof}

Let a map $M \rightarrow K$ act on a butterfly as the maps $Y \rightarrow M \rightarrow K \rightarrow X'$ and $Y' \rightarrow M \rightarrow K \rightarrow X$ compatibly act on the two extensions in the above product.

Consider an automorphism of a butterfly \eqref{diamondbutterfly}. Subtracting the identity yields a map between the entire diagram which is zero except for a map $Q \rightarrow Q$. Each such map must factor uniquely as $Q \rightarrow M \rightarrow K \rightarrow Q$. This shows that 2-isomorphisms between two fixed butterflies in $\IIsom(\xi, \eta)$ are a pseudo-torsor under $\hh{K}$.

We show local existence. Localizing in $M$, we assume $Y \simeq M \oplus P$ and $Y' \simeq M \oplus P'$ are split extensions. Define $Q := (X' \oplus_K X) \oplus M$, with butterfly diagram: 
\begin{equation}\label{locexistsbutterfly}
\begin{tikzcd}
        &       &P \arrow{dr}     &       &       \\
        &X \arrow{ur} \arrow{dr}     &       &M \oplus P \arrow{dr}     &       \\
K \arrow{ur} \arrow{dr}       &       &Q \arrow{ur} \arrow{dr}      &       &M      \\
        &X' \arrow{ur} \arrow{dr}     &       &M \oplus P' \arrow{ur}     &       \\
        &       &P' \arrow{ur}     &       &       
\end{tikzcd}\end{equation} 

In order to show butterflies are pairwise locally isomorphic, pick an arbitrary butterfly $Q'$ filling in the above diagram. Localizing in $M$ sufficiently, $Q'$ splits as $(X \oplus_K X') \oplus M$; we may choose an isomorphism of $Q'$ with the above $Q$ compatible with all the structure maps. 

We leave the verification that all relevant composites in \eqref{locexistsbutterfly} are short exact sequences and that the diagrams formed by our map of butterflies commute to the dedicated reader. 

\end{proof}

We can finally describe the isomorphism of Theorem \ref{extiscohom} in the case $p = 2$.

\begin{proposition}\label{explicitsplittingsp2}
The map $\Split : \Ext_A^2(M, K) \simeq H^2(\Amod/M, \hh{K})$ furnished by Theorem \ref{extiscohom} sends a 2-extension to its $\hh{K}$-gerbe of splittings. 
\end{proposition}

\begin{proof}

Write 
\[m :   \quad      0 \rightarrow      P \rightarrow      Y \rightarrow      M \rightarrow        0      \]
\[\gamma : \quad   0 \rightarrow     K \rightarrow      X \rightarrow      P \rightarrow      0  \]
so that $\xi  = \gamma \smile m$

Consider the long exact sequence in $\Ext^p$ and $H^p$ coming from $\gamma$. The isomorphism of Theorem \ref{extiscohom} is one of universal $\delta$-functors, so we get a commutative diagram: 
\[\begin{tikzcd}
\cdots \arrow{r}     &\Ext_A^1(M, X) \arrow{r} \arrow[phantom, "\circ"]{dr} \arrow{d}{\sim}        &\Ext_A^1(M, P) \arrow{r}{\gamma \smile} \arrow{d}{\sim}  \arrow[phantom, "\circ"]{dr}      &\Ext_A^2(M, K) \arrow{d}{\sim}     \\
\cdots \arrow{r}      &H^1(M, \hh{X}) \arrow{r}     &H^1(M, \hh{P}) \arrow{r}     &H^2(M, \hh{K})
\end{tikzcd}\]

The extension $m$ maps to the 2-extension $\xi$ under the boundary map $\gamma \smile$ by definition. To see what $\xi$ maps to in $H^2$, send $m$ around the bottom corner of the square. The boundary map on cohomology sends the torsor $\hh{Y}$ associated to $m$ to its $\hh{K}$-gerbe of lifts to an $\hh{X}$-torsor. By commutativity of the left square, this is equivalent to the $\hh{K}$-gerbe of lifts of the extension $m$ to an extension by $X$. As depicted in the rearranged butterfly diagram below, this gerbe is identical to $\Split(\xi)$. 

\begin{equation}
\begin{tikzcd}
        &       &P \arrow{dr}       &       &       \\
        &X \arrow{ur} \arrow[dashed]{dr}       &       &Y \arrow{dr}       &       \\
K \arrow[equals]{dr} \arrow{ur}        &       &Q \arrow[dashed]{ur} \arrow[dashed]{dr}       &       &M       \\
        &K \arrow[dashed]{ur} \arrow{dr}       &       &M \arrow[equals]{ur}       &       \\
        &       &0 \arrow{ur}       &       &       
\end{tikzcd}
\end{equation}

\end{proof}

\hspace{10pc}

The final ingredient in Theorem \ref{obstcomparison} is the map 
\[\Hom_A(J \otimes M, K) \overset{\smile \omega}{\rightarrow} \Ext_A^2(M, K)\] 

\noindent of Diagram \eqref{obstcomparisondiagram} and \eqref{illusieexactseq}. This homomorphism sends $u : J \otimes M \rightarrow K$ to its pushout $u \smile \omega$ along a fixed 2-extension $\omega$. 

In order to construct $\omega$, take a flat $A'$-module mapping surjectively onto $M$ with kernel $L$: 
\[0 \rightarrow L \rightarrow H \rightarrow M \rightarrow 0.\]

When we tensor with $A$, we obtain 

\[\omega : \quad 0 \rightarrow J \otimes M \rightarrow \overline{L} \rightarrow \overline{H} \rightarrow M \rightarrow 0.\]

We write $\overline{L}$ for $L \otimes A = L/JL$. Since $\omega$ computes $\Tor_0^{A'}(A, M)$ and $\Tor_1^{A'}(A, M)$, it's the canonical obstruction $\omega(A', M)$ in Illusie's work by \cite[IV.3.1.9]{illusie1}.

We must check $\omega$ is well-defined up to isomorphism. Given two flat surjections onto $M$, we can always choose a third surjecting onto both (e.g., the direct sum of a cover by free $A'$-modules trivializing both extensions). We may assume there is a map between the two flat resolutions: 
\[\begin{tikzcd}
0 \arrow{r}      &L'\arrow{r} \arrow{d}    &H' \arrow{r} \arrow{d}   &M \arrow{r} \arrow[equals]{d}    &0      \\
0 \arrow{r}      &L \arrow{r}     &H \arrow{r}     &M \arrow{r}     &0
\end{tikzcd}\]

In this case, simply tensor the whole diagram by $A$ to get a map of complexes between the two definitions of $\omega$. A morphism between 2-extensions as chain complexes induces a butterfly as in Paragraph \ref{inducedbutterfly} in the appendix. 

Hence $\omega$ is sufficiently well-defined to define a morphism to the group of connected components $\Ext_A^2(M, K)$, even though there's no canonical complex-level representative. The reader is free to fix one representative $\omega$ and transpose to a given one via the above. 

\hspace{10pc}

\openrefcount{obstcomparison}
\begin{theorem}\label{newobstcomparison}
The diagram 
\begin{equation}
\begin{tikzcd}
\Hom_A(J \otimes M, K) \arrow{dr}[swap]{\DDef} \arrow{r}{\smile \omega}      &\Ext_A^2(M, K) \arrow{d}{\Split}     \\
        &H^2(\Amod/M, \hh{K})
\end{tikzcd}
\end{equation}
anti-commutes. 
\end{theorem}
\closerefcount

Fix $u \in \Hom_A(J \otimes M, K)$ and continue to write $\overline{L} := L \otimes A$. Our proof consists of two lemmas: one exhibits a functor $\beta : \Split(u \smile \omega) \rightarrow \DDef(\inpt, u, K)$ over $\Amod/M$, and the other shows it's an $\hh{K}$-anti-equivalence.

\begin{lemma}\label{constructafunctor}
There is a natural functor $\beta : \Split(u \smile \omega) \rightarrow \DDef(\inpt, u, K)$ over $\Amod/M$. 
\end{lemma}

\begin{proof}

Given a splitting (Remark \ref{splitpushoutbutterfly}): 

\begin{equation}\label{splittingbutterflydeformation}
\begin{tikzcd}[column sep=small, row sep=small]
0 \arrow{rr}       &&J \otimes M \arrow{rr} \arrow{dd}{u}      &&\overline{L} \arrow{rr} \arrow{dr}       &&\overline{H} \arrow{rr}       &&M \arrow{rr} \arrow[equals]{dd}      &&0     \\
        &&      &&       &Q \arrow{dr} \arrow{ur}       &      &&     &&\\
0 \arrow{rr}       &&K \arrow[equals]{rr}     &&K \arrow{rr}{0} \arrow{ur}     &&M \arrow[equals]{rr}     &&M \arrow{rr}     &&0
\end{tikzcd}\end{equation}

of $u \smile \omega$, consider the pushout
\[\begin{tikzcd}
\eta : \arrow{dd}     &0 \arrow{r}       &L \arrow{d} \arrow{r} \arrow[phantom, very near end, "\lrcorner"]{ddr}      &H \arrow{r} \arrow{dd}      &M \arrow{r} \arrow[equals]{dd}      &0      \\
        &&\overline{L} \arrow{d}       &&&\\
\zeta :     &0 \arrow{r}       &Q \arrow{r}      &H \oplus_L Q \arrow{r}{h}       &M \arrow{r}      &0
\end{tikzcd}\]

Distinguish between three natural maps $H \oplus_L Q \rightarrow M$: 

\begin{itemize}

\item $H \oplus_L Q \overset{h}{\rightarrow} M$ is the structure map $H \rightarrow M$ and zero on $Q$. 

\item $H \oplus_L Q \overset{q}{\rightarrow} M$ is the structure map $Q \rightarrow M$ and zero on $H$. 

\item $H \oplus_L Q \overset{h + q}{\rightarrow} M$ is the sum of the two maps above, given by both structure maps. It factors through $\overline{H}$.
\end{itemize}

Define an extension $\xi$ by taking the kernel 
\begin{equation}\label{definem'xizeta}
\begin{tikzcd}
\xi : \arrow{d}      &0 \arrow{r}      &K \arrow[dashed]{r} \arrow{d}      &M' \arrow[dashed, "g"]{r} \arrow[dashed]{d}     &M \arrow{r} \arrow[equals]{d}     &0      \\
\zeta :     &0 \arrow{r}      &Q \arrow{r} \arrow{d}      &H \oplus_L Q \arrow{d}{h+q} \arrow{r}{h}       &M \arrow{r}     &0      \\
&0 \arrow{r}       &\overline{H} \arrow[equals]{r}       &\overline{H} \arrow{r}       &0
\end{tikzcd}\end{equation}

To show $\theta(\xi) = u$, tensor the diagram $\eta \rightarrow \zeta \leftarrow \xi$ by $\otimes_{A'} A$. We get a diagram containing 
\[\begin{tikzcd}
\eta \otimes A : \arrow{d}      &0 \arrow{r}      &J \otimes M \arrow{r} \arrow[equals]{d}      &\overline{L} \arrow{r} \arrow{d}       &\cdots         \\
\zeta \otimes A :     &\cdots \arrow{r}      &J \otimes M \arrow{r} \arrow[equals]{d}      &Q \arrow{r}       &\cdots         \\
\xi \otimes A : \arrow{u}       &\cdots \arrow{r}      &J \otimes M \arrow{r}{\theta(\xi)}      &K \arrow{r} \arrow[hook]{u}      &\cdots         
\end{tikzcd}\]

The left pentagon of the original butterfly verifies that the map $J \otimes M \rightarrow \overline{L} \rightarrow Q$ factors as $J \otimes M \overset{u}{\rightarrow} K \hookrightarrow Q$. Combine this with the above diagram into 

\[\begin{tikzcd}
        &J \otimes M \arrow{dl}[swap]{\theta(\xi)} \arrow{d} \arrow{dr}{u}        &       \\
K \arrow{dr}       &\overline{L} \arrow{d}       &K \arrow{dl}      \\
        &Q
\end{tikzcd}\]

Hence $J \otimes M \overset{u}{\rightarrow} K \rightarrow Q$ and $J \otimes M \overset{\theta(\xi)}{\rightarrow} K \rightarrow Q$ are the same. Since $K \hookrightarrow Q$ is a monomorphism, this confirms $u = \theta(\xi)$.

An isomorphism $Q \simeq Q'$ of butterflies induces a unique isomorphism $H \oplus_L Q \simeq H \oplus_L Q'$ fixing $H$ and $L$. These isomorphisms are compatible with a functorial isomorphism of the whole diagram \eqref{definem'xizeta} inducing the identity on $K$, $M$, and $\overline{H}$, whence a unique isomorphism on kernels $M' \simeq M''$. Let this be the action of $\beta$ on arrows. 

\end{proof}

\begin{lemma}
The functor $\beta$ of Lemma \ref{constructafunctor} is an anti-equivalence. 
\end{lemma}

\begin{proof}

We continue to use terminology from the proof of Lemma \ref{constructafunctor}. 

A morphism of gerbes which is banded by an isomorphism is an equivalence by \cite[IV.2.2.7]{giraudnonab}. We claim that $\beta$ is banded not by the identity, but by $-id_K$. 

A map $M \overset{\varphi}{\rightarrow} K$ acts on a butterfly \eqref{splittingbutterflydeformation} by adding the map $Q \rightarrow \overline{H} \rightarrow M \overset{\varphi}{\rightarrow} K \rightarrow Q$ to the identity on $Q$. Then the induced automorphism of $H \oplus_L Q$ is obtained by adding the identity to 
\begin{equation}\label{hplusqautom} H \oplus_L Q \overset{q}{\rightarrow} M \overset{\varphi}{\rightarrow} K \rightarrow H \oplus_L Q\end{equation}

We claim ``a," ``b," and ``c" in the following solid diagram commute: 

\begin{equation}\label{actiononextensions}
\begin{tikzcd}
0 \arrow{r}       &M' \arrow{r} \arrow[swap, "-g"]{dr} \arrow[bend left, phantom, "a"]{dr} \arrow[dashed, swap, "z"]{ddd}      &H \oplus_L Q \arrow{d}{q} \arrow[phantom, "c"]{dddr} \arrow{r}       &\overline{H} \arrow{ddd}{0} \arrow{r}       &0      \\
        &       &M \arrow{d}{\varphi}      &       &       \\
        &       &K \arrow{dl} \arrow{d} \arrow[bend left, phantom, "b"]{dl}      &       &       \\
0 \arrow{r}       &M' \arrow{r}     &H \oplus_L Q \arrow{r}       &\overline{H} \arrow{r}       &0
\end{tikzcd}\end{equation}

The map $-g$ is the additive inverse of $g : M' \rightarrow M$ defined in \eqref{definem'xizeta}, and $z$ is the composite. 

Diagram \eqref{definem'xizeta} witnesses the commutativity of rectangle ``c" and triangle ``b." The same diagram also asserts $M' \rightarrow H \oplus_L Q \overset{h+q}{\rightarrow} \overline{H} \rightarrow M$ is zero; equivalently, that 
\[-g : M' \rightarrow H \oplus_L Q \overset{-h}{\rightarrow} M\]
and 
\[M' \rightarrow H \oplus_L Q \overset{q}{\rightarrow} M\]
are equal. This confirms commutativity of triangle ``a."

Diagram \eqref{actiononextensions} defines a morphism of complexes; add the identity morphism to obtain a morphism of complexes given by 
\begin{itemize}
\item $id$ on $\overline{H}$. 
\item $id + $\eqref{hplusqautom} on $H \oplus_L Q$. 
\item $id + z$ on $M'$. 
\end{itemize}

By unwinding the definition of $\beta$ on arrows, the action of $\varphi$ on $Q$ is sent by $\beta$ to the automorphism of $\xi$ which is $id + z$ on $M'$ and the identity on $M$ and $K$. This is precisely the action of $-\varphi$ on $\xi$. 

\end{proof}

%\section{Deformations and Butterflies}

%\subfile{4zbutterflies}

\section{Illusie's Exact Sequence}\label{illusiecomparison}

%Show that Illusie's exact sequence is exact and comes from the Grothendieck-Leray spectral sequence

This section describes Illusie's Exact Sequence

\begin{equation}
0 \rightarrow \Ext_A^1(M, K) \rightarrow \Ext_{A'}^1(M, K) \overset{\theta}{\rightarrow} \Hom_A(J \otimes M, K) \overset{\smile \omega}{\rightarrow} \Ext_A^2(M, K) \tag{\ref{illusieexactseq}}
\end{equation}

\noindent and proves Theorem \ref{illusie5term}. 

We need naturality to construct the comparison diagram \eqref{illusie5termiso}. 

\begin{lemma}\label{illusieexactseqisnatural}
The maps in \eqref{illusieexactseq} are all natural in $M$ and $K$. 
\end{lemma}

\begin{proof}

Naturality of the first arrow in \eqref{illusieexactseq} is clear.

Consider the pushout and pullback of an extension $\xi \in \Ext_{A'}^1(M, K)$ along maps $K \rightarrow L$ and $N \rightarrow M$. 

\[\begin{tikzcd}
0 \arrow{r}      &K \arrow{r} \arrow[equals]{d}    &M' \times_M N \arrow{r} \arrow[very near start, phantom, "\ulcorner"]{dr} \arrow{d}     &N \arrow{r}  \arrow{d}   &0      \\
0 \arrow{r}      &K \arrow[phantom, very near end, "\lrcorner"]{dr}\arrow{r} \arrow{d}    &M' \arrow{r} \arrow{d}     &M \arrow{r}  \arrow[equals]{d}   &0      \\
0 \arrow{r}      &L \arrow{r}     &M' \oplus_K L \arrow{r}     &M \arrow{r}     &0      
\end{tikzcd}\]

Tensoring $\inpt \otimes A$, we get

\[\begin{tikzcd}
J \otimes N \arrow{d} \arrow{r}     &K \arrow[equals]{d}  \arrow{r}    &\cdots\\
J \otimes M \arrow[equals]{d} \arrow{r}{\theta(\xi)}     &K \arrow{d} \arrow{r}     &\cdots\\
J \otimes M \arrow{r}     &L \arrow{r}     &\cdots
\end{tikzcd}\]

Commutativity of this diagram implies that $\theta(\xi|_N) = \theta(\xi) \circ (J \otimes N \rightarrow J \otimes M)$ and $\theta((K \rightarrow L) \smile \xi) = (K \rightarrow L) \circ \theta(\xi)$. Conclude the arrow $\theta$ in Diagram \eqref{illusieexactseq} is natural.

If $\theta(\xi) = 0$, or equivalently if $M'$ is an $A$-module, then $\theta(\xi|_N) = \theta((K \rightarrow L) \smile \xi) = 0$ and $M' \times_M N$ and $M' \oplus_K L$ are both $A$-modules. The inclusion $\Ext_A^1(M, K) \subseteq \Ext_{A'}^1(M, K)$ beginning the sequence is natural. 

The associativity of pushing out and pulling back 2-extensions furnishes the naturality of the last arrow, $\smile \omega$.

\end{proof}

\begin{remark}
Let $A + \epsilon M$ be the trivial squarezero algebra extension of $A$ by $M$: the $A$-module $A \oplus M$ endowed with multiplication given by $A$'s action and $M$ squaring to zero. It may be graded by placing $M$ in degree 1, $A$ in degree 0. 

Illusie defined the exact sequence \eqref{illusieexactseq} using the first graded piece of the cotangent complex transitivity triangle $\mathbb{L}^{\text{gr}}_{A+\epsilon M/A/A'}$. The compatibility of that approach with this more direct one was verified already by Illusie as follows: 
$$\begin{tikzcd}
0 \arrow{r}     &\Ext_A^1(M, K) \arrow{r}   \arrow[equals]{d}    &\Ext_A^1(\mathbb{L}^{\text{gr}}_{A+\epsilon M/A'}, K) \arrow{r} \arrow{d}{\sim}     &\Ext_A^1(\mathbb{L}_{A/A'}\otimes_A M, K) \arrow{r} \arrow{d}{\sim}      &\Ext_A^2(M, K)  \arrow[equals]{d}     \\
0 \arrow{r}   &\Ext_A^1(M, K) \arrow{r}             &\Ext_{A'}^1(M, K) \arrow{r}       &\Hom_{A}(J \otimes M, K) \arrow{r}{\smile \omega}       &\Ext_A^2(M, K) 
\end{tikzcd}$$
The commutativity of the leftmost square was observed immediately before Proposition 3.1.5 \cite[pg. 248]{illusie1}, and the middle square is equivalent to diagram (3.1.3) on the previous page. 

In the rightmost square, we are cupping with $\omega$. Its two middle terms are $\tau_{[-1}(M \overset{L}{\otimes} A)$ by construction. By the naturality of both sequences in $K$ (the top is obtained by applying the functor $\Ext_A^1(\inpt, K)$ to the transitivity triangle), it suffices to reduce to the case where $K := J \otimes M$ and $u  = id_{J \otimes M}$. Since $\omega$ is the ``canonical obstruction'' by IV.3.1.9 of \cite[pg. 250]{illusie1}, the square commutes. (For us, the ground ring $\Upsilon$ is $A'$.)

\end{remark}

It remains to show the following diagram commutes, and to describe the dashed arrow. 

\begin{equation}
\begin{tikzcd}
0 \arrow{r}       &\Ext_A^1(M, K) \arrow{r} \arrow{d}{\sim}     &\Ext_{A'}^1(M, K) \arrow{r} \arrow{d}{\sim}    &\Hom_A(J \otimes M, K) \arrow{r} \arrow[dashed, "\sim"]{d}    &\Ext_A^2(M, K) \arrow{d}{\Split}             \\
0 \arrow{r}       &H^1(A/M, \hh{K}) \arrow{r}      &H^1(A'/M, \hh{K}) \arrow{r}     &H^0(A/M, R^1\pi_* \hh{K}) \arrow{r}       &H^2(A/M, \hh{K}) 
\end{tikzcd}\tag{\ref{illusie5termiso}}
\end{equation}

\begin{lemma}\label{illusie5termisnatural}
The solid arrows in Diagram \eqref{illusie5termiso} are natural in the $A$-modules $M$ and $K$. 
\end{lemma}

\begin{proof}

The whole solid diagram is natural in $K$ by Lemmas \ref{illusieexactseqisnatural}, \ref{grothssinterpretation}, and by construction of the isomorphism provided by Theorem \ref{extiscohom}. The same is true of $M$, except possibly the vertical isomorphisms. Choose $N \rightarrow M$ in $\Amod$. 

Given an extension 
\[0 \rightarrow K \rightarrow M' \rightarrow M \rightarrow 0\]

\noindent representing the torsor $\hh{K} \action \hh{M'|M}$, the pulled back extension represents $\hh{M' \times_M N | N} \simeq j^*\hh{M'}$. The naturality in $H^1$ comes from pullback of torsors. This trick shows the naturality of the first two vertical isomorphisms.

Now examine a 2-extension 
\[\xi : \quad       0 \rightarrow K \rightarrow X \rightarrow Y \rightarrow M \rightarrow 0\]

A section of $\Split(\xi|_N)$ over $T \rightarrow N$ is simply a section of $\Split(\xi)$ over $T \rightarrow N \rightarrow M$, so we have a (strict) 2-fiber product
\[\begin{tikzcd}
\Split(\xi|_N) \arrow[phantom, very near start, "\ulcorner"]{dr}\arrow{r} \arrow{d}      &\Split(\xi) \arrow{d}        \\
\Amod/N \arrow{r}{j_!}     &\Amod/M
\end{tikzcd}\]

Accordingly, $\Split(\xi|_N) \simeq \Split(\xi)|_N$, where the first pullback belongs to $\Ext_A^2$ and the second to $H^2$. Hence $\Split$ is natural in $M$.

\end{proof}

\begin{lemma}
The dashed arrow in Diagram \eqref{illusie5termiso} exists and is an isomorphism. 
\end{lemma}

\begin{proof}

Lemma \ref{illusie5termisnatural} shows that the diagram is natural in $M \in \Amod$. Sheafify to obtain 

\[\begin{tikzcd}
\Ext_{A'}^1(\inpt, K)^{sh} \arrow{r} \arrow{d}         &\Hom_A(J \otimes \inpt, K) \arrow[dashed]{d}     \\
H^1(A'/\inpt, \hh{K})^{sh} \arrow{r}       &R^1 \pi_* \hh{K}
\end{tikzcd}\]

All of the arrows are isomorphisms, since the outer terms go to zero. Define the sought-after isomorphism as the composition of the other three.

\end{proof}

\begin{remark}
This last argument describes the maps $\Ext_{A'}^1(M, K) \overset{\theta}{\rightarrow} \Hom_A(J \otimes M, K)$ and $H^1(A'/M, \hh{K}) \rightarrow H^0(A/M, R^1 \pi_* \hh{K})$ as sheafification.
\end{remark}

\openrefcount{illusie5term}
\begin{theorem}
Diagram \eqref{illusie5termiso} commutes. Scilicet, Illusie's exact sequence \eqref{illusieexactseq} and the 5-term exact sequence from the Grothendieck-Leray spectral sequence are isomorphic. 
\end{theorem}
\closerefcount

\begin{proof}

Take an extension $\xi : 0 \rightarrow K \rightarrow M' \rightarrow M \rightarrow 0$ of $A$-modules. Whether one first considers $K$, $M'$, and $M$ as $A'$-modules and then forms the torsor $\hh{K} \action \hh{M'}$ on $\AAmod/M$ or forms the torsor $\hh{K} \action \hh{M'}$ on $\Amod/M$ and then applies $\pi^*$ makes no difference: $\pi^*\hh{K} = \hh{K}$ functorially. The left square commutes. 

To verify that the two maps to a sheaf $\Ext_{A'}^1(\inpt, K) \rightrightarrows R^1 \pi_* \hh{K}$ agree, we may sheafify. Then the arrow $\Hom_A(J \otimes \inpt, K) \rightarrow R^1 \pi_* \hh{K}$ was defined to make this square commute.

The rightmost square remains. We show

\[\begin{tikzcd}
\Hom_A(J \otimes M, K) \arrow{r} \arrow{d} \arrow{dr}{-\DDef}      &\Ext_A^2(M, K) \arrow{d}{\Split}     \\
H^0(A/M, R^1 \pi_* \hh{K}) \arrow{r}      &H^2(A/M, \hh{K})
\end{tikzcd}\]
commutes. The upper right triangle commutes by Theorem \ref{obstcomparison}. 

Under the lower left triangle, consider the image of $u \in \Hom_A(J \otimes M, K)$ under the two maps. The bottom horizontal arrow sends a global section to the inverse of (the class of) its gerbe of lifts to an $\hh{K}$-torsor in $\AAmod/M$ by \ref{grothssinterpretation} in the appendix. 

By the commutativity of the leftmost square of diagram \eqref{illusie5termiso}, we see that this gerbe is equivalent to the gerbe of lifts of the corresponding map $\Hom_A(J \otimes M, K)$ to an $A'$-module extension. This was the definition of $\DDef(\inpt, u, K)$.

\end{proof}

%\section{Conventions and Notation}

%\subfile{6conventions}

%Appendix with technical details and maybe the vanishing lemma
\appendix

\section{The Grothendieck Spectral Sequence}

In this appendix, we describe a map belonging to the 5-term exact sequence induced from the Grothendieck Spectral Sequence in the special case of a morphism of topoi. 

We also specify our sign convention for torsors and gerbes: 

\begin{definition}\label{torsorgerbesignconvention}
Given a short exact sequence
\[0 \rightarrow A \rightarrow B \overset{g}{\rightarrow} C \rightarrow 0\]
\noindent of abelian groups in $E$, consider the corresponding long exact sequence in sheaf cohomology: 
\[\begin{tikzcd}
\cdots \arrow{r}      &H^0(B) \arrow{r}{g_*}     &H^0(C) \arrow{r}{\partial^0}     &H^1(A) \arrow{r}     &\cdots\\
\cdots \arrow{r}      &H^1(B) \arrow{r}{g_*}     &H^1(C) \arrow{r}{\partial^1}     &H^2(A) \arrow{r}     &\cdots
\end{tikzcd}\]

We agree that 
\begin{itemize}
\item The boundary map $\partial^0$ sends $\gamma \in H^0(C)$ to the $A$-torsor $P$ whose sections over $U \in E$ are
\[P(U) := \{\beta \in H^0(U, B) \:\:|\:\: g_*\beta = \gamma|_U\}\]
\item The boundary map $\partial^1$ sends the class of a $C$-torsor $P$ to the $A$-gerbe whose sections over $U \in E$ are $B|_U$-torsors $Q$ on $E/U$ with a $g$-equivariant map $Q \rightarrow P|_U$. The arrows are $B$-maps.  
\end{itemize}

\end{definition}

Another convention is to choose the boundary maps $-\partial^i$ instead. 

Now consider $\pi : X \rightarrow Y$ and $\rho : Y \rightarrow Z$ two functors between abelian categories with enough injectives and suppose $\pi$ sends injective objects to $\rho$-acyclic ones. Given $G \in X$, we construct the Grothendieck Spectral Sequence as follows (see the diagrams below): 
\begin{itemize}

\item Resolve $G$ by injectives, $\{J^p\}$. 

\item Apply $\pi$. 

\item Find a Cartan-Eilenberg resolution $\{I^{p, q}\}$ of the resulting complex and apply $\rho$ to it. 

\item Take horizontal ($p$) and then vertical ($q$) cohomology to get $E_2^{p, q} = R^q \rho R^p \pi G$. 

\end{itemize}

Remark 13.21.4 of \cite[015G]{sta} observes the naturality of Cartan-Eilenberg resolutions with respect to maps of chain complexes $J^\bullet$, and any choice of injective resolutions is natural in $G$. The Grothendieck Spectral Sequence is then natural in $G$ via the usual functoriality of the pair of right-derived functors. 

We use the convention of \cite[012X]{sta} that the Total Complex should have differentials $d_\rightarrow^{p,q} + (-1)^p d_\uparrow^{p,q}$, where the horizontal differential $d_\rightarrow$ has degree $(+1, 0)$ and the vertical $d_\uparrow$ degree $(0, +1)$.

Our case of interest is when $X$ and $Y$ are ringed topoi, $\pi = f_*$ for a morphism $f$ of topoi, and $\rho = \Gamma$. We construct this spectral sequence.  

\[\begin{tikzcd}
I^{03}      &       &       &       \\
I^{02} \arrow{u} \arrow{r}    &I^{12}     &       &       \\
I^{01} \arrow{u} \arrow{r}    &I^{11} \arrow{u}{d_\uparrow} \arrow{r}{d_\rightarrow}    &I^{21}     &       \\
I^{00} \arrow{u} \arrow{r}      &I^{10} \arrow{u} \arrow{r}    &I^{20}  \arrow{u} \arrow{r}   &I^{30}      
\end{tikzcd}
E_1 : 
\left[\begin{tikzcd}
h_\rightarrow^{03}      &       &       &       \\
h_\rightarrow^{02} \arrow{u}     &h_\rightarrow^{12}     &       &       \\
h_\rightarrow^{01} \arrow{u}     &h_\rightarrow^{11} \arrow{u}{d_\uparrow}    &h_\rightarrow^{21}     &       \\
h_\rightarrow^{00} \arrow{u}       &h_\rightarrow^{10} \arrow{u}     &h_\rightarrow^{20}  \arrow{u}    &h_\rightarrow^{30}      
\end{tikzcd}\right.\]

\[E_2 : 
\left[\begin{tikzcd}
H^3(f_*(G))      &       &       &       \\
H^2(f_*(G))      &H^2(R^1f_*(G))     &       &       \\
H^1(f_*(G))      &H^1(R^1f_*(G)) \arrow{ulu}{d_2}    &H^1(R^2f_*(G))     &       \\
\Gamma(f_*(G))        &\Gamma(R^1f_*(G)) \arrow{ulu}     &\Gamma(R^2f_*(G))  \arrow{uul}    &\Gamma(R^3f_*(G))      
\end{tikzcd}\right.\]

A section of the derived pushforward $\Gamma(R^1 f_* G)$ may be thought of as an element $\beta \in \Gamma(I^{10})$ such that $d_\rightarrow \beta = 0$ and $d_\uparrow \beta$ is in the image of $d_\rightarrow$. Choose a lift $\gamma \in \Gamma(I^{01})$ of $d_\uparrow \beta$. 

Suppose another section $\gamma'$ maps to $d_\uparrow \beta$. Then $\gamma-\gamma'$ lies in $h^{01}_\rightarrow$. Since $H^2(f_* G)$ is a quotient by the image of $h^{01}_\rightarrow$, the image $d_\uparrow \gamma = d_\uparrow \gamma'$ is a well-defined cohomology class.

Consider the two subsheaves of $I^{01}$: 

\[P(U) := \{s \in \Gamma(U, I^{01}) \: | \: d_\rightarrow s = d_\uparrow \beta \text{ and } d_\uparrow s = 0\}\]
and
\[Q(U) := \{t \in \Gamma(U, I^{01}) \: | \: d_\uparrow t = d_\uparrow \gamma \text{ and } d_\rightarrow t = 0\}\]

Let $L$ denote the kernel of both differentials $d_\uparrow$ and $d_\rightarrow$ in $I^{01}$. Then $P$ and $Q$ are naturally $L$-torsors. The sum of elements from $P$ and $Q$ is the $L$-torsor of local sections which map to $d_\rightarrow \beta$ horizontally and $d_\uparrow \gamma$ vertically; $\gamma$ trivializes this torsor. 

\[\begin{tikzcd}
L \oplus L \arrow{r}{+} \arrow[symbol=\circlearrowleft]{d}      &L \arrow[symbol=\circlearrowleft]{d}      \\
P \oplus Q \arrow{r}{+} \arrow[hook]{d}      &L \arrow[hook]{d}      \\
I^{01} \oplus I^{01} \arrow{r}{+}        &I^{01}
\end{tikzcd}\]

This diagram witnesses that the sum of the torsors $P$ and $Q$ is zero, hence that they are inverses. Since we have an exact sequence 
\[0 \rightarrow f_*G \rightarrow I^{00} \rightarrow L \rightarrow 0\]
with middle term injective, we may identify $L$-torsors with $f_*G$-gerbes via the boundary map. 

Remark that, by definition, $P$ was the $f_*G$-gerbe of lifts of $\beta$ to a $G$-torsor on $Y$. Likewise, $Q$ was the image of $\beta$ under the $E_2$-page differential. That is, the map $\Gamma(R^1f_*G) \rightarrow H^2(f_* G)$ in the spectral sequence and the corresponding 5-term exact sequence sends a global section to the inverse of the gerbe of its lifts to a torsor on $Y$.

We package the observations made in this section into a lemma. 

\begin{lemma}\label{grothssinterpretation}

The Grothendieck Spectral Sequence is natural in $G$. The map $\Gamma(R^1 f_*G) \rightarrow H^2(f_*G)$ sends a global section to the inverse of the class of its gerbe of lifts to a $G$-torsor on $Y$.

\end{lemma}

%%%%%%%%%%%%%%%%%%%%Butterflies

\section{Butterflies}\label{butterfliesappendix}

Several essential properties of the 2-groupoid of 2-extensions are transcribed from \cite{butterflies} or \cite{sga7i} to our context. In the process, we fix many of the conventions requisite for working with them. 

\begin{definition}\label{butterfly2categorydefinitions}

We define and begin our study of $\EExt_A^2(\inpt, K) \rightarrow \Amod$.

The 2-groupoid $\EExt_A^2(\inpt, K) \rightarrow \Amod/M$ has: 

\begin{itemize}
\item Sections: 2-extensions 
\[0 \rightarrow K \rightarrow X \rightarrow Y \rightarrow T \rightarrow 0\] 
\noindent over $T \rightarrow M$. 

\item Morphisms: Butterflies inducing the identity on $K$: 
\begin{equation}\label{butterflydiagram}
\begin{tikzcd}
0 \arrow{r}      &K \arrow{r} \arrow[equals]{dd}     &X\arrow{dr} \arrow{rr}     &      &Y \arrow{r}       &T  \arrow{dd}\arrow{r}    &0      \\
        &       &       &Q \arrow{dr} \arrow{ur}      &       &       &       \\
0 \arrow{r}      &K \arrow{r}     &X'\arrow{rr} \arrow{ur}     &      &Y' \arrow{r}      &T'  \arrow{r}    &0
\end{tikzcd}\end{equation}
\noindent lie over the morphism $T \rightarrow T'$. 

\item 2-Isomorphisms: Maps between the extensions defining the two butterflies. 
\end{itemize}
\end{definition}

\paragraph{Induced Butterfly}\label{inducedbutterfly}
Given a map of chain complexes that are 2-extensions: 

$$\begin{tikzcd}
0 \arrow{r}      &K \arrow{r} \arrow[equals]{d}     &X\arrow{d} \arrow{r}           &Y \arrow{r} \arrow{d}       &T  \arrow{d}\arrow{r}    &0      \\
0 \arrow{r}      &K \arrow{r}     &X'\arrow{r}           &Y' \arrow{r}      &T'  \arrow{r}    &0
\end{tikzcd}$$

\noindent we get a butterfly: 

$$\begin{tikzcd}
0 \arrow{r}      &K \arrow{r}\arrow[equals]{dd}     &X \arrow{rr} \arrow{dr}      &&Y \arrow{r}     &T \arrow{r}\arrow{dd}     &0      \\
        &       &       &X' \oplus Y \arrow{dr} \arrow{ur}       &       &       &\\
0 \arrow{r}      &K \arrow{r}     &X' \arrow{rr} \arrow{ur}     &&Y' \arrow{r}     &T' \arrow{r}     &0
\end{tikzcd}$$

The map $X \rightarrow X' \oplus Y$ is the sum of the two maps $X \rightarrow X'$ and $X \rightarrow Y$; the map $X' \oplus Y \rightarrow Y'$ is $X' \rightarrow Y'$ \textit{minus} $Y \rightarrow Y'$. Define the identity to be the induced butterfly from the identity map on chain complexes.

\paragraph{Composition of Butterflies}
Define composition as follows: 

$$\begin{tikzcd}
0 \arrow{r}      &K \arrow{r} \arrow[equals]{dd}     &X\arrow{dr} \arrow{rr}     &      &Y \arrow{r}       &T  \arrow{dd}\arrow{r}    &0      \\
        &       &       &Q \arrow{dr} \arrow{ur}      &       &       &       \\
0 \arrow{r}      &K \arrow{r} \arrow[equals]{dd}     &X'\arrow{dr} \arrow{ur} \arrow{rr}     &      &Y' \arrow{r}       &T'  \arrow{dd}\arrow{r}    &0      \\
        &       &       &Q' \arrow{dr} \arrow{ur}      &       &       &       \\
0 \arrow{r}      &K \arrow{r}     &X''\arrow{rr} \arrow{ur}     &      &Y'' \arrow{r}      &T''  \arrow{r}    &0
\end{tikzcd}$$

\noindent compose to the butterfly: 
\[\begin{tikzcd}
0 \arrow{r}      &K \arrow{r} \arrow[equals]{dd}     &X\arrow{dr} \arrow{rr}     &      &Y \arrow{r}       &T  \arrow{dd}\arrow{r}    &0      \\
        &       &       &Q\oplus_{Y'}^{X'}Q' \arrow{dr} \arrow{ur}      &       &       &       \\
0 \arrow{r}      &K \arrow{r}     &X''\arrow{rr} \arrow{ur}     &      &Y'' \arrow{r}      &T''  \arrow{r}    &0
\end{tikzcd}\]

Here, the term $Q \oplus_{Y'}^{X'} Q'$ refers to the cokernel of $X' \rightarrow Q \times_{Y'} Q'$. The SW-NE diagonal sequence is the cokernel of the map of extensions
\[\begin{tikzcd}
        &0  \arrow{r}     &X'\arrow[equals]{r} \arrow{d}    &X'\arrow{r} \arrow{d}    &0      \\
0 \arrow{r}      &X'' \arrow{r}       &Q \times_{Y'} Q'\arrow{r}       &Q \arrow{r}     &0
\end{tikzcd}\]

It is exact by the snake lemma. Exactness of the NW-SE diagonal follows analogously from the following lemma.

\begin{lemma}
The module $Q \oplus_{Y'}^{X'} Q'$ is the kernel of $Q \amalg_{X'} Q' \rightarrow Y'$. 
\end{lemma}

\begin{proof}
Consider the pair of diagrams: 

\[\begin{tikzcd}
0 \arrow{r}   &X'' \arrow{r} \arrow[equals]{d}    &Q \times_{Y'} Q' \arrow{r} \arrow{d} \arrow[phantom, very near start, "\ulcorner"]{dr}      &Q\arrow{r} \arrow{d}     &0      \\
0 \arrow{r}  &X''\arrow{r}     &Q'\arrow{r}     &Y' \arrow{r}       &0
\end{tikzcd}\quad
\begin{tikzcd}
0 \arrow{r}   &X' \arrow{r} \arrow{d}    &Q \arrow{r} \arrow{d}       &Y\arrow{r} \arrow[equals]{d}     &0      \\
0 \arrow{r}  &Q'\arrow{r}     &Q \amalg_{X'} Q'\arrow{r}  \arrow[phantom, very near start, "\lrcorner"]{ul}   &Y \arrow{r}       &0
\end{tikzcd}.\]

Let $\pi_i$ denote the projections of $Q \times_{Y'} Q'$ and $\iota_i$ denote the inclusions of $Q \amalg_{X'} Q'$, $i = 1,2$. We get a map from the pullback to the pushout above: 
\[\begin{tikzcd}
0 \arrow{r}      &X'' \arrow{r} \arrow{d} \arrow[phantom, near start, "\circ"]{dr}    &Q \times_{Y'} Q' \arrow{dl}{\pi_2} \arrow{r}{\pi_1} \arrow[phantom, near end, "\circ"]{dr} &Q \arrow{r} \arrow{d} \arrow{dl}[swap]{\iota_1}    &0      \\
0 \arrow{r}      &Q' \arrow{r}[swap]{-\iota_2}    &Q \amalg_{X'} Q' \arrow{r}      &Y \arrow{r}     &0
\end{tikzcd}\]

The left and right vertical maps are structure maps coming from the butterfly. The left and right triangles commute because they are the left and right squares of the pullback and pushout diagrams above, respectively. We don't claim the inner parallelogram is commutative. 

The arrows $\pi_2$ and $\iota_1$ define a nullhomotopy of a morphism of chain complexes. We place the nullhomotopic morphism $\psi = \iota_1 \circ \pi_1 - \iota_2 \circ \pi_2$ in the middle of a diagram: 
\[\begin{tikzcd}
       &0 \arrow{r}      &X' \arrow[equals]{r} \arrow{d}      &X' \arrow{r} \arrow{d}       &0       \\
0 \arrow{r}      &X'' \arrow{r} \arrow{d}     &Q \times_{Y'} Q'  \arrow{r} \arrow{d}{\psi}  &Q \arrow{r} \arrow{d}     &0      \\
0 \arrow{r}      &Q' \arrow{r} \arrow{d}    &Q \amalg_{X'} Q' \arrow{r} \arrow{d}      &Y \arrow{r}     &0      \\
0 \arrow{r}       &Y' \arrow[equals]{r}     &Y' \arrow{r}     &0
\end{tikzcd}\]

Then either of the vertical composites of chain maps in the above diagram is zero, and we get a map from the cokernel of the first chain map to the kernel of the last: 

\[\begin{tikzcd}
0 \arrow{r}       &X'' \arrow[equals]{d} \arrow{r}     &Q \oplus_{Y'}^{X'} Q' \arrow{r} \arrow{d}      &Y \arrow{r} \arrow[equals]{d}     &0       \\
0 \arrow{r}      &X''  \arrow{r}      &\ker \arrow{r}      &Y\arrow{r}      &0
\end{tikzcd}\]

Here $\ker = \ker (Q \amalg_{X'} Q' \rightarrow Y')$. The careful reader will notice that the structure map $X'' \rightarrow \ker$ is the negative of the usual map because we used $-\iota_2$, but we still have an isomorphism $Q \oplus_{Y'}^{X'} Q' \simeq \ker$. 

\end{proof}

\begin{lemma}
Suppose $T \rightarrow T'$ is an isomorphism in Diagram \eqref{butterflydiagram}. The butterfly given by flipping the diagram upside-down is its inverse (up to 2-isomorphism). 
\end{lemma}

\begin{proof}
The composite is 

\[\begin{tikzcd}
0 \arrow{r}      &K \arrow{r} \arrow[equals]{dd}     &X\arrow{dr} \arrow{rr}     &      &Y \arrow{r}       &T  \arrow{d}\arrow{r}    &0      \\
        &       &       &Q\oplus_{Y'}^{X'}Q \arrow{dr} \arrow{ur}      &       &T' \arrow{d}       &       \\
0 \arrow{r}      &K \arrow{r}     &X \arrow{rr} \arrow{ur}     &      &Y \arrow{r}      &T  \arrow{r}    &0
\end{tikzcd}\]

The diagonal $Q \rightarrow Q \times_{Y'} Q$ gives a section of the middle row of the diagram constructing $Q \oplus_{Y'}^{X'} Q$:  
\[\begin{tikzcd}
        &0 \arrow{r}     &X' \arrow[equals]{r} \arrow{d}    &X'\arrow{r} \arrow{d}    &0      \\
0 \arrow{r}      &X \arrow{r} \arrow{d}    &Q \times_{Y'} Q \arrow{r} \arrow{d}      &Q \arrow{r} \arrow{d} \arrow[bend right, dashed]{l} \arrow[dashed]{dl}    &0      \\
0 \arrow{r}      &X \arrow{r}     &Q \oplus_{Y'}^{X'} Q \arrow{r}        &Y \arrow{r} \arrow[bend right, dashed]{l}     &0
\end{tikzcd}\]

Compose to get the map $Q \dashrightarrow Q \oplus_{Y'}^{X'} Q$. The map $X' \rightarrow Q \times_{Y'} Q$ factors through the diagonal, since the maps $X' \rightarrow Q$ are the same. 

Then the map $X' \rightarrow Q \dashrightarrow Q \oplus_{Y'}^{X'} Q$ factors as $X' \rightarrow Q \times_{Y'} Q \rightarrow Q \oplus_{Y'}^{X'} Q$, which is the zero map. Hence $Q \dashrightarrow Q \oplus_{Y'}^{X'} Q$ factors through the cokernel, inducing a section $Y \dashrightarrow Q \oplus_{Y'}^{X'} Q$. 

The exact sequence defining the composite butterfly splits, so it's isomorphic to the identity butterfly with $X \oplus Y$ in the center. 

\end{proof}

\begin{lemma}

Up to 2-isomorphism, butterflies $\xi \simeq \eta$ over $N \rightarrow M \in \Amod$ are the same as butterflies $\xi \simeq \eta|_N$ over $id_N$.

\end{lemma}

\begin{proof}

Let $\xi|_N$ denote the 2-extension 
\[0 \rightarrow K \rightarrow X \rightarrow Y|_N \rightarrow N \rightarrow 0\]
and so forth.

Given a butterfly 

\[\begin{tikzcd}
0 \arrow{r}      &K \arrow{r} \arrow[equals]{dd}     &X\arrow{dr} \arrow{rr}     &      &Y \arrow{r}       &N  \arrow{dd}\arrow{r}    &0      \\
        &       &       &Q \arrow{dr} \arrow{ur}      &       &        &       \\
0 \arrow{r}      &K \arrow{r}     &X' \arrow{rr} \arrow{ur}     &      &Y' \arrow{r}      &M  \arrow{r}    &0
\end{tikzcd}\]

\noindent over $N \rightarrow M \in \Amod$, 

\[\begin{tikzcd}
K \arrow{r} \arrow{d}       &X \arrow{r} \arrow{d}      &P \arrow{d}     \\
X' \arrow{r} \arrow{d}     &Q \arrow{r} \arrow[dashed]{d}     &Y \arrow{d}     \\
P' \arrow{r} \arrow[equals]{d}     &Y'|_N \arrow{r} \arrow{d} \arrow[phantom, very near start, "\ulcorner"]{dr}     &N  \arrow{d}    \\
P' \arrow{r}     &Y' \arrow{r}     &M
\end{tikzcd}\]

The dashed arrow comes from the rest of the diagram. Discard the bottom row to get a 3x3 grid of modules whose rows are known to be exact. The 3x3 Lemma ensures that the middle column is exact. We can rearrange to obtain a butterfly $\xi \simeq \eta|_N$. 

The composite butterfly $\xi \simeq \eta|_N \simeq \eta$ is then 

\[\begin{tikzcd}
0 \arrow{r}      &K \arrow{r} \arrow[equals]{dd}     &X\arrow{dr} \arrow{rr}     &      &Y \arrow{r}       &N  \arrow{dd}\arrow{r}    &0      \\
        &       &       &(Q \oplus X')/X' \arrow{dr} \arrow{ur}      &       &       &       \\
0 \arrow{r}      &K \arrow{r}     &X' \arrow{rr} \arrow{ur}     &      &Y' \arrow{r}      &M  \arrow{r}    &0
\end{tikzcd}\]

\noindent where $(Q \oplus X')/X'$ includes $X'$ via the sum of the two natural maps. Unwinding definitions carefully, we see that the diagonal arrows are: 
\begin{itemize}
\item $X' \rightarrow (Q \oplus X')/X'$ includes the second summand. 
\item $X \rightarrow (Q \oplus X')/X'$ is the structure map $X \rightarrow Q$. 
\item $Q \oplus X' \rightarrow Y'|_N \oplus X' \rightarrow Y'$, where the first arrow is the pair of natural maps and the second is $X' \rightarrow Y'$ \textit{minus} $Y'|_N \rightarrow Y'$. 
\item $Q \oplus X' \rightarrow Q \rightarrow Y$. 
\end{itemize}
Precomposing the latter two maps by the sum of the natural maps $X' \rightarrow Q \oplus X'$ gives zero, yielding the factorization through $(Q \oplus X')/X'$.

The reader may check commutativity of the following diagrams.  

\[\begin{tikzcd}[column sep=small, row sep=small]
        &(Q \oplus X')/X' \arrow{dd} \arrow{dr}       &       \\
X' \arrow{ur} \arrow{dr}      &       &Y      \\
        &Q \arrow{ur}
\end{tikzcd}
\begin{tikzcd}[column sep=small, row sep=small]
        &(Q \oplus X')/X' \arrow{dd} \arrow{dr}       &       \\
X \arrow{ur} \arrow{dr}      &       &Y'      \\
        &Q \arrow{ur}
\end{tikzcd}\]

The map $Q \oplus X' \rightarrow Q$ is the map $X' \rightarrow Q$ \textit{minus} $id_Q$. The 3-Lemma says this arrow is an isomorphism, and the two diagrams together build a 2-isomorphism between the two butterflies.

\end{proof}

\begin{remark}\label{splitpushoutbutterfly}

By a symmetric proof, butterflies $\xi \simeq \eta$

\[\begin{tikzcd}
\xi:        &0 \arrow{r}      &K \arrow{r} \arrow{dd}     &X\arrow{dr} \arrow{rr}     &      &Y \arrow{r}       &M  \arrow[equals]{dd}\arrow{r}    &0      \\
&        &       &       &Q \arrow{dr} \arrow{ur}      &       &        &       \\
\eta :      &0 \arrow{r}      &L \arrow{r}     &X' \arrow{rr} \arrow{ur}     &      &Y' \arrow{r}      &M  \arrow{r}    &0
\end{tikzcd}\]

\noindent which allow non-identity left vertical arrows are the same as butterflies $(K \rightarrow L) \smile \xi \simeq \eta$ with left vertical arrow $id_L$.

\end{remark}

We provide a reassurance that the group $\Ext_A^2(M, K)$ of connected components of $\EExt_A^2(M, K)$ is the same, whether the maps are butterflies or maps of chain complexes. 

\begin{lemma}

The connected components of $\EExt_A^2(M, K)$ under butterflies are the same as the connected components of 2-extensions under morphisms of complexes. 

\end{lemma}

\begin{proof}

If we have a morphism of complexes between two 2-extensions, we get an induced butterfly as above. Conversely, suppose we have a butterfly between two 2-extensions. Choose a cover $\{N_i \rightarrow M\}$ which trivializes the butterfly, available due to Lemma \ref{isomisagerbe}. Let $T := \bigoplus N_i$. Then the butterfly pulls back under $T \rightarrow M$ to the split butterfly. Butterflies come from morphisms of chain complexes if and only if they're split, so we have morphisms of chain complexes connecting the two 2-extensions which shared a butterfly. 

\[\begin{tikzcd}
0  \arrow{r}     &K \arrow{r}       &X \arrow{r}      &Y \arrow{r}      &M  \arrow{r}     &0      \\
0 \arrow{r}      &K \arrow{r} \arrow[equals]{d} \arrow[equals]{u}    &X \arrow{r} \arrow[equals]{u} \arrow{d}     &Y|_T  \arrow{u}\arrow{r} \arrow{d}     &T  \arrow{u}\arrow{r} \arrow[equals]{d}     &0      \\
0  \arrow{r}     &K \arrow{r} \arrow[equals]{d}     &X'\arrow{r} \arrow[equals]{d}     &Y'|_T  \arrow{r} \arrow{d}    &T \arrow{r} \arrow{d}     &0      \\
0 \arrow{r}      &K \arrow{r}     &X' \arrow{r}    &Y' \arrow{r}    &M \arrow{r}     &0
\end{tikzcd}\]

\end{proof}

%\section{DELETE ME}

%\subfile{excerptcechextandcohom}

\bibliographystyle{unsrt}%Used BibTeX style is unsrt
\bibliography{bib}

\end{document}